\documentclass[11pt]{amsart}
\usepackage[T1]{fontenc}
\usepackage[utf8]{inputenc}
\usepackage{lmodern}
\usepackage{amsmath,amssymb,amsthm,mathtools}
\usepackage{mathrsfs}
\usepackage{enumitem}
\usepackage{microtype}
\usepackage[margin=1.08in]{geometry}
\usepackage{amsaddr}
\usepackage[
colorlinks=true,
linkcolor=blue,
citecolor=blue,
urlcolor=blue
]{hyperref}
\usepackage{bookmark}

\numberwithin{equation}{section}

\newtheorem{theorem}{Theorem}[section]
\newtheorem{proposition}[theorem]{Proposition}
\newtheorem{lemma}[theorem]{Lemma}

\newtheorem{corollary}[theorem]{Corollary}
\theoremstyle{definition}

\theoremstyle{remark}

\newcommand{\R}{\mathbb R}

\newcommand{\intt}{\operatorname{int}}

\newcommand{\spanop}{\operatorname{span}}
\newcommand{\diag}{\operatorname{diag}}
\newcommand{\Id}{\mathrm{Id}}
\newcommand{\norm}[1]{\lVert#1\rVert}
\newcommand{\ip}[2]{\langle#1,#2\rangle}

\newcommand{\ma}{\mathrm{ma}}

\newcommand{\one}{\mathbf1}

\setlist[enumerate]{
	label=\textup{(\roman*)},
	leftmargin=*
}
\setlist[itemize]{leftmargin=*}

\makeatletter
\renewcommand{\@set@authors@addresses}{

	\@setauthors

	\begingroup
	\centering
	\def\author##1{}
	\def\\{\protect\linebreak}
	\def\address##1##2{
		\par
		{\normalfont\small\itshape\ignorespaces##2\par}
	}
	\def\email##1##2{}
	\def\curraddr##1##2{}
	\def\urladdr##1##2{}
	\addresses
	\par
	\endgroup
}
\makeatother

\author{Yuxian Dong}
\author{Wenle Shi\textsuperscript{*}}
\address{School of Mathematical Sciences,
	Dalian University of Technology,\\
Dalian 116024, P. R. China}

\thanks{\textsuperscript{*}Corresponding author.}
\thanks{
	\textit{E-mail addresses:}
	\nolinkurl{yuxiand@hotmail.com} (Y. Dong);
	\nolinkurl{shi-wenle@hotmail.com} (W. Shi).
}



\title{Some polynomial spaces are Euclidean balls}

\begin{document}
\begin{abstract}
In this paper, we establish a matrix criterion for a space to be
homeomorphic to a closed Euclidean ball.
As applications, we prove that the normalized space of
$K$-Lorentzian polynomials on an arbitrary open convex cone $K$
is homeomorphic to a closed Euclidean ball, giving an affirmative
answer to Br\"and\'en and Leake's question, and we also settle Br\"and\'en's conjecture that the normalized space of homogeneous multiaffine real stable polynomials is homeomorphic to a closed Euclidean ball.
In addition, our results show that the normalized space of
homogeneous $K$-stable polynomials of degree $d$ is homeomorphic
to a closed Euclidean ball.
\end{abstract}
\maketitle

\section{Introduction}\label{sec:intro}

Let $n\ge1$ be an integer and write $\mathbf{x}=(x_1,\ldots,x_n)$.
Let $\mathbb R[\mathbf{x}]$ denote the ring of real polynomials in the variables $x_1,\ldots,x_n$.
For an integer $d\ge1$, a polynomial $f\in\mathbb R[\mathbf{x}]$ is called \emph{homogeneous of degree $d$} if
\[
f(t\mathbf{x})=t^d f(\mathbf{x}) \qquad (\mathbf{x}\in\mathbb R^n,\,t\in\mathbb R).
\]
Let $\mathrm H_d(\mathbb R^n)$ denote the vector space of real homogeneous polynomials of degree $d$, and set $\mathrm H_0(\mathbb R^n)=\mathbb R$.

For each $d\ge0$, fix an ordering $m_1,\ldots,m_M$ of the monomials of total degree $d$ in $x_1,\ldots,x_n$, where $M=\binom{n+d-1}{d}$.
Define the coefficient map
\[
\operatorname{coef}_d:\mathrm H_d(\mathbb R^n)\longrightarrow\mathbb R^M,
\qquad
\operatorname{coef}_d\!\left(\sum_{j=1}^{M}a_jm_j\right)=(a_1,\ldots,a_M).
\]
This map is a linear isomorphism.
For $f,g\in\mathrm H_d(\mathbb R^n)$, write $\operatorname{coef}_d(f)=(a_1,\ldots,a_M)$ and $\operatorname{coef}_d(g)=(b_1,\ldots,b_M)$.
Define a metric $\rho$ on $\mathrm H_d(\mathbb R^n)$ by
\[
\rho(f,g)=\left(\sum_{j=1}^{M}(a_j-b_j)^2 \right)^{1/2}.
\]
We equip $\mathrm H_d(\mathbb R^n)$ with the topology induced by the metric $\rho$.
For an integer $m\ge1$, the \emph{closed Euclidean unit ball}
in $\mathbb R^m$ is
\[
\overline B^m
=
\left\{
\mathbf{y}=(y_1,\ldots,y_m)\in\mathbb R^m:
\sum_{i=1}^{m}y_i^2\le1
\right\}.
\]
Let $\overline B^0$ denote a one-point space.
Let $X$ and $Y$ be topological spaces.
For a map $h:X\to Y$ and a subset $U\subseteq Y$, the \emph{inverse image} of $U$ under $h$ is
\[
h^{-1}(U)=\{x\in X:h(x)\in U\}.
\]
The map $h$ is called \emph{continuous} if $h^{-1}(U)$ is open in $X$ for every open subset $U$ of $Y$.
Let $X$ and $Y$ be topological spaces, and let $h:X\to Y$ be a bijection.
If both $h$ and its inverse $h^{-1}:Y\to X$ are continuous, then $h$ is called a \emph{homeomorphism}.
The spaces $X$ and $Y$ are said to be \emph{homeomorphic} if there exists a homeomorphism from $X$ to $Y$.

A nonempty open set $K\subseteq\mathbb R^n$ is called an \emph{open convex cone} if
\[
a\mathbf{u}+b\mathbf{v}\in K
\quad\text{for all }\mathbf{u},\mathbf{v}\in K
\text{ and }a,b>0.
\]
We assume throughout that $K\ne\mathbb R^n$.

For $\mathbf{u}=(u_1,\ldots,u_n)\in\mathbb R^n$ and $f\in\mathbb R[\mathbf{x}]$, the \emph{directional derivative} of $f$ in the direction $\mathbf{u}$ is defined by
\[
D_{\mathbf{u}}f(\mathbf{x})
=
\left.\frac{d}{dt}f(\mathbf{x}+t\mathbf{u})\right|_{t=0}
=
\sum_{i=1}^{n}u_i\frac{\partial f}{\partial x_i}(\mathbf{x}).
\]
For $f\in\mathbb R[\mathbf{x}]$, its \emph{Hessian matrix} is defined by
\[
\operatorname{Hess}(f)
=
\left(\frac{\partial^2f}{\partial x_i\partial x_j}\right)_{i,j=1}^{n}.
\]
When $f\in\mathrm H_2(\mathbb R^n)$, this is a constant real symmetric matrix.

Let $d\ge2$.
Following Br\"and\'en and Leake~\cite[Definition~2.1]{BL}, a nonzero polynomial $f\in\mathrm H_d(\mathbb R^n)$ is called \emph{$K$-Lorentzian} if the following two conditions hold:
\begin{enumerate}
\item[(i)]
$D_{\mathbf{u}_1}\cdots D_{\mathbf{u}_d}f>0$ for all $\mathbf{u}_1,\ldots,\mathbf{u}_d\in K$.
\item[(ii)]
For all $\mathbf{u}_3,\ldots,\mathbf{u}_d\in K$, the matrix $\operatorname{Hess}\bigl(D_{\mathbf{u}_3}\cdots D_{\mathbf{u}_d}f\bigr)$ has exactly one positive eigenvalue.
\end{enumerate}
When $d=2$, the matrix in condition~(ii) is $\operatorname{Hess}(f)$.
For $d=1$, a polynomial $f\in\mathrm H_1(\mathbb R^n)$ is called $K$-Lorentzian if $f(\mathbf{u})>0$ for every $\mathbf{u}\in K$.
For $d=0$, a constant polynomial $f\equiv c$ is called $K$-Lorentzian if $c\ge0$.

Fix $\mathbf{v}=(v_1,\ldots,v_n)\in K$. 
For $d\ge0$, define
\[
\mathcal L_{K,d,\mathbf{v}}
=
\left\{
f\in\mathrm H_d(\mathbb R^n):
f\text{ is }K\text{-Lorentzian},\ f(\mathbf{v})=1
\right\}.
\]
Let $\mathbb R_{>0}^n=\{\mathbf{x}\in\mathbb R^n: x_i>0\text{ for }1\le i\le n\}$ and let $\mathbf{1}=(1,\ldots,1)\in\mathbb R^n$.
By~\cite[Proposition~8.2]{BL}, $\mathcal L_{\mathbb R_{>0}^n,d,\mathbf{1}}$ is the space of Lorentzian polynomials $f\in\mathrm H_d(\mathbb R^n)$ satisfying $f(\mathbf{1})=1$.

Br\"and\'en and Huh~\cite{BH} proved that $\mathcal L_{\mathbb R_{>0}^n,d,\mathbf{1}}$ is compact and contractible.
They further conjectured that this space is homeomorphic to a closed Euclidean ball.
Br\"and\'en~\cite{BrandenBalls} subsequently proved this conjecture.
For an arbitrary open convex cone $K$, Br\"and\'en and
Leake~\cite[Section~9]{BL} proved that $\mathcal L_{K,d,\mathbf v}$ is compact and contractible.
They asked whether this space is homeomorphic to a closed Euclidean ball~\cite[Question~1]{BL}.

Let $\overline K$ denote the closure of $K$.
The \emph{lineality space} of $\overline K$ is $S=\overline K\cap(-\overline K)$.
Set $r=n-\dim S$.

\begin{theorem}[{\cite[Question~1]{BL}}]
\label{thm:lor-main}
For every integer $d\ge0$, the space $\mathcal L_{K,d,\mathbf v}$ is homeomorphic to a closed Euclidean ball of dimension $\binom{r+d-1}{d}-1$.
\end{theorem}

We next consider $K$-stable polynomials.
A nonzero polynomial $f\in\mathbb R[\mathbf{x}]$ is called \emph{$K$-stable} if $f(\mathbf{x}+\mathrm{i}\mathbf{y})\ne0$ for all $\mathbf{x}\in\mathbb R^n$ and $\mathbf{y}\in K$, where $\mathrm{i}^2=-1$.
For $d\ge0$, define
\[
\mathcal S_{K,d,\mathbf{v}}
=
\left\{
f\in\mathrm H_d(\mathbb R^n):
f\text{ is }K\text{-stable},\ f(\mathbf{v})=1
\right\}.
\]
A polynomial $f\in\mathbb R[\mathbf{x}]$ is called \emph{real stable} if it is $\mathbb R_{>0}^n$-stable.
Br\"and\'en~\cite[Theorem~4.3]{BrandenBalls} proved that $\mathcal S_{\mathbb R_{>0}^n,d,\mathbf{1}}$ is homeomorphic to a closed Euclidean ball.
The following theorem extends this result to arbitrary open convex cones.
\begin{theorem}
\label{thm:stable-main}
For every integer $d\ge0$, the space $\mathcal S_{K,d,\mathbf v}$ is homeomorphic to a closed Euclidean ball of dimension $\binom{r+d-1}{d}-1$.
\end{theorem}

A polynomial $f\in\mathbb R[\mathbf{x}]$ is called \emph{multiaffine} if its degree in each variable is at most one.
For an integer $d$ with $0\le d\le n$, define
\[
\mathcal M_{n,d}
=
\left\{
f\in\mathcal S_{\mathbb R_{>0}^n,d,\mathbf{1}}:
f\text{ is multiaffine}
\right\}.
\]
Br\"and\'en conjectured that $\mathcal M_{n,d}$ is homeomorphic to a closed Euclidean ball.

\begin{theorem}[{\cite[Conjecture~5.5]{BrandenBalls}}]
\label{thm:ma-main}
For every integer $d$ with $0\le d\le n$, the space $\mathcal M_{n,d}$ is homeomorphic to a closed Euclidean ball of dimension $\binom{n}{d}-1$.
\end{theorem}

We prove Theorems~\ref{thm:lor-main}, \ref{thm:stable-main},
and~\ref{thm:ma-main} using the following matrix criterion.
We first introduce the necessary topological notation.
Let $m\ge1$ be an integer.
A set $U\subseteq\mathbb R^m$ is called \emph{closed} if $\mathbb R^m\setminus U$ is open.
The \emph{interior} and \emph{closure} of $U$ are defined by
\begin{align*}
\operatorname{int}U
&=
\bigcup\{O\subseteq U:
O\text{ is open in }\mathbb R^m\},\\
\overline U
&=
\bigcap\{F\subseteq\mathbb R^m:
U\subseteq F,\ F\text{ is closed in }\mathbb R^m\}.
\end{align*}
The \emph{boundary} of $U$ is $\partial U=\overline U\setminus\operatorname{int}U$.
For $R>0$, write
\[
R\overline B^m
=
\left\{
\mathbf y=(y_1,\ldots,y_m)\in\mathbb R^m:
\sum_{i=1}^m y_i^2\le R^2
\right\}.
\]
The set $U$ is called \emph{bounded} if $U\subseteq R\overline B^m$ for some $R>0$.
Let $S^{m-1}=\partial\overline B^m$ denote the Euclidean unit sphere.

\begin{theorem}
\label{thm:main-criterion}
Let $m\ge1$ be an integer and $C,Y\subseteq\mathbb R^m$ be closed sets such that $C$ is bounded and $X=C\cap Y$ is nonempty.
Fix $D=\operatorname{diag}(d_1,\ldots,d_m)$, $B\in\mathbb R^{m\times m}$, and $\mathbf c\in\mathbb R^m$, where $d_i>0$ for $1\le i\le m$.
Suppose that for every $T>0$ and every continuous function $\lambda:[0,T]\to[0,\infty)$, every solution $\mathbf u:[0,T]\to\mathbb R^m$ of the linear differential equation
\[
\mathbf u'(t)
=
-D\mathbf u(t)
+\lambda(t)\bigl(B\mathbf u(t)+\mathbf c\bigr)
\]
satisfies the following conditions:
\begin{enumerate}
\item If $\lambda\equiv0$ and $\mathbf u(0)\in C$, then $\mathbf u(T)\in\operatorname{int}C$.
\item If $\lambda(t)>0$ for every $t\in[0,T]$ and $\mathbf u(0)\in Y$, then $\mathbf u(T)\in\operatorname{int}Y$.
\end{enumerate}
Then there exists a homeomorphism $h:X\to\overline B^m$ such that $h(\partial X)=S^{m-1}$.
\end{theorem}

The paper is organized as follows.
Section~\ref{sec:flows} proves Theorem~\ref{thm:main-criterion}.
Section~\ref{sec:cones} proves Theorems~\ref{thm:lor-main} and~\ref{thm:stable-main} for Lorentzian and stable polynomials on cones.
Section~\ref{sec:ma} proves Theorem~\ref{thm:ma-main} for multiaffine real stable polynomials.
The linear algebra needed for the multiaffine case is developed in Appendix~\ref{app:matrix}.

\section{A matrix criterion for closed balls}
\label{sec:flows}

In this section, we prove Theorem~\ref{thm:main-criterion} using the contractive flow criterion of Galashin, Karp, and Lam~\cite{GKL}.
Let $m\ge1$ be an integer.
For $\mathbf x=(x_1,\ldots,x_m)$ and $\mathbf y=(y_1,\ldots,y_m)$ in $\mathbb R^m$, the Euclidean inner product and norm are defined by
\[
\langle\mathbf x,\mathbf y\rangle
=
\sum_{i=1}^{m}x_i y_i,
\qquad
\|\mathbf x\|
=
\left(\sum_{i=1}^{m}x_i^2\right)^{1/2}.
\]
For a matrix $A\in\mathbb R^{m\times m}$, the \emph{operator norm} induced by the Euclidean norm is defined by
\[
\|A\|_{\mathrm{op}}
=
\sup\{\|A\mathbf z\|:
\mathbf z\in\mathbb R^m,\ \|\mathbf z\|=1\}.
\]
For every $\mathbf z\in\mathbb R^m$, we have $\|A\mathbf z\|\le\|A\|_{\mathrm{op}}\|\mathbf z\|$.

A \emph{flow} on $\mathbb R^m$ is a continuous map $F:\mathbb R\times\mathbb R^m\to\mathbb R^m$ satisfying
\[
F(0,\mathbf z)=\mathbf z,
\qquad
F(s+t,\mathbf z)=F(s,F(t,\mathbf z))
\]
for all $s,t\in\mathbb R$ and $\mathbf z\in\mathbb R^m$.
Following Galashin, Karp, and Lam~\cite[Definition~2.1]{GKL}, we call a flow $F$ \emph{contractive} if
\[
\|F(t,\mathbf z)\|<\|\mathbf z\|
\qquad(t>0,\,\mathbf z\in\mathbb R^m\setminus\{\mathbf 0\}).
\]

\begin{lemma}[{\cite[Lemma~2.3]{GKL}}]
\label{lem:gkl}
Let $U\subseteq\mathbb R^m$ be a nonempty bounded open set, and let $F:\mathbb R\times\mathbb R^m\to\mathbb R^m$ be a contractive flow.
Suppose that
\[
F(t,\mathbf z)\in U
\qquad(t>0,\, \mathbf z\in\overline U).
\]
Then there exists a homeomorphism $h:\overline U\to\overline B^m$ such that $h(\partial U)=S^{m-1}$.
\end{lemma}

A \emph{vector field} on $\mathbb R^m$ is a map $b:\mathbb R^m\to\mathbb R^m$.
We say that $b$ is \emph{Lipschitz} if there exists a constant $L\ge0$ such that
\[
\|b(\mathbf x)-b(\mathbf y)\|\le L\|\mathbf x-\mathbf y\|
\qquad(\mathbf x,\mathbf y\in\mathbb R^m).
\]

Suppose that $b$ is Lipschitz.
By~\cite[Theorem~2.2 and Corollary~2.6]{Teschl}, for each $\mathbf z\in\mathbb R^m$ the initial value problem
\[
\mathbf u'(t)=b(\mathbf u(t)),\qquad \mathbf u(0)=\mathbf z
\]
has a unique solution $\mathbf u_{\mathbf z}:\mathbb R\to\mathbb R^m$.
Define $\Phi:\mathbb R\times\mathbb R^m\to\mathbb R^m$ by $\Phi(t,\mathbf z)=\mathbf u_{\mathbf z}(t)$.
The map $\Phi$ is continuous by~\cite[(2.43) and Theorem~2.9]{Teschl}.
The initial condition gives $\Phi(0,\mathbf z)=\mathbf z$, and uniqueness of solutions to the autonomous equation gives $\Phi(s+t,\mathbf z)=\Phi(s,\Phi(t,\mathbf z))$ for all $s,t\in\mathbb R$ and $\mathbf z\in\mathbb R^m$.
We call $\Phi$ the \emph{flow generated by $b$}.

\begin{lemma}\label{lem:boundary-robustness}
Let $C\subseteq\mathbb R^m$ be closed, and let $b_0$ and $b$ be Lipschitz vector fields on $\mathbb R^m$ with flows $T$ and $\Phi$, respectively.
Suppose that $b=b_0$ on $\mathbb R^m\setminus\operatorname{int}C$ and that $T(t,\mathbf x)\in\operatorname{int}C$ for all $t>0$ and $\mathbf x\in C$.
Then $\Phi(t,\mathbf x)\in\operatorname{int}C$ for all $t>0$ and $\mathbf x\in C$.
\end{lemma}

\begin{proof}
Fix $\mathbf y\in\mathbb R^m\setminus\operatorname{int}C$.
We claim that $T(s,\mathbf y)\notin C$ for all $s<0$.
Suppose, to the contrary, that $T(s,\mathbf y)\in C$ for some $s<0$.
Since $-s>0$, we have
\[
\mathbf y=T(-s,T(s,\mathbf y))\in\operatorname{int}C,
\]
contrary to the choice of $\mathbf y$.

For all $s\le0$, we have $T(s,\mathbf y)\notin\operatorname{int}C$.
Since $b=b_0$ on $\mathbb R^m\setminus\operatorname{int}C$, we obtain
\[
\frac{d}{ds}T(s,\mathbf y)
=
b_0(T(s,\mathbf y))
=
b(T(s,\mathbf y))
\qquad(s\le0).
\]
Both $T(0,\mathbf y)$ and $\Phi(0,\mathbf y)$ equal $\mathbf y$.
By uniqueness of solutions,
$\Phi(s,\mathbf y)=T(s,\mathbf y)$ for all $s\le0$.
In particular, $\Phi(-t,\mathbf y)\notin C$ for every $t>0$.

Now fix $\mathbf x\in C$ and $t>0$.
Suppose, to the contrary, that $\Phi(t,\mathbf x)\notin\operatorname{int}C$.
By the preceding conclusion, $\Phi(-t,\Phi(t,\mathbf x))\notin C$.
However, the flow identity gives
\[
\Phi(-t,\Phi(t,\mathbf x))
=
\Phi(0,\mathbf x)
=
\mathbf x\in C,
\]
a contradiction.
\end{proof}

\begin{proof}[Proof of Theorem~\ref{thm:main-criterion}]
Since $C$ is bounded and $m\ge1$, the complement $\mathbb R^m\setminus C$ is nonempty.
Define $\delta_C:\mathbb R^m\to[0,\infty)$ by
\begin{equation}\label{eq:distance-cutoff}
\delta_C(\mathbf x)
=
\inf_{\mathbf z\in\mathbb R^m\setminus C}
\|\mathbf x-\mathbf z\|.
\end{equation}
For $\mathbf x,\mathbf y\in\mathbb R^m$ and $\mathbf z\notin C$, by the triangle inequality, we have
\[
\delta_C(\mathbf x)
\le \|\mathbf x-\mathbf z\|
\le \|\mathbf x-\mathbf y\|+\|\mathbf y-\mathbf z\|.
\]
Taking the infimum over $\mathbf z\notin C$, we obtain
\[
|\delta_C(\mathbf x)-\delta_C(\mathbf y)|
\le \|\mathbf x-\mathbf y\|.
\]
Hence $\delta_C$ is Lipschitz with constant one.
If $\mathbf x\in\operatorname{int}C$, there exists $\eta>0$ such that $\{\mathbf z\in\mathbb R^m:\|\mathbf z-\mathbf x\|<\eta\}
\subseteq C$.
By the definition of $\delta_C$, we have $\delta_C(\mathbf x)\ge\eta>0$.
If $\mathbf x\notin\operatorname{int}C$, every open ball centered at $\mathbf x$ meets $\mathbb R^m\setminus C$, and therefore $\delta_C(\mathbf x)=0$.

Define $G:\mathbb R^m\to\mathbb R^m$ by
\[
G(\mathbf x)=\delta_C(\mathbf x)(B\mathbf x+\mathbf c).
\]
Since $C$ is nonempty and bounded and $\delta_C$ is Lipschitz, the number
\[
L_G
=
\sup_{\mathbf z\in C}\|B\mathbf z+\mathbf c\|
+
\|B\|_{\mathrm{op}}\sup_{\mathbf z\in C}\delta_C(\mathbf z)
\]
is finite.
For $\mathbf x\in C$ and $\mathbf y\in\mathbb R^m$, we have
\[
G(\mathbf x)-G(\mathbf y)
=
\bigl(\delta_C(\mathbf x)-\delta_C(\mathbf y)\bigr)
(B\mathbf x+\mathbf c)+\delta_C(\mathbf y)B(\mathbf x-\mathbf y).
\]
Since $\delta_C$ vanishes outside $C$, we have
$0\le\delta_C(\mathbf y)\le\sup_{\mathbf z\in C}\delta_C(\mathbf z)$
for every $\mathbf y\in\mathbb R^m$.
By the Lipschitz estimate for $\delta_C$, we have
\[
\|G(\mathbf x)-G(\mathbf y)\|
\le
\bigl(\|B\mathbf x+\mathbf c\|
+\delta_C(\mathbf y)\|B\|_{\mathrm{op}}\bigr)
\|\mathbf x-\mathbf y\|\le L_G\|\mathbf x-\mathbf y\|.
\]
By symmetry, the same estimate holds when $\mathbf y\in C$.
If both $\mathbf x$ and $\mathbf y$ lie outside $C$, then
$G(\mathbf x)=G(\mathbf y)=\mathbf0$.
Hence $G$ is Lipschitz on $\mathbb R^m$ with constant $L_G$.

Choose $\varepsilon>0$ such that
\[
\kappa=\min_{1\le i\le m}d_i-\varepsilon L_G>0,
\]
and define $b:\mathbb R^m\to\mathbb R^m$ by
\begin{equation}\label{eq:matrix-global-field}
b(\mathbf x)=-D\mathbf x+\varepsilon G(\mathbf x).
\end{equation}
The vector field $b$ is Lipschitz.
Let $\Phi:\mathbb R\times\mathbb R^m\to\mathbb R^m$ be the flow generated by $b$.
By the diagonal form of $D$ and the Lipschitz estimate for $G$, for all $\mathbf x,\mathbf y\in\mathbb R^m$, we have
\[
\begin{aligned}
\langle b(\mathbf x)-b(\mathbf y),\mathbf x-\mathbf y\rangle
&=
-\langle D(\mathbf x-\mathbf y),\mathbf x-\mathbf y\rangle\\
&\quad+\varepsilon\langle G(\mathbf x)-G(\mathbf y),
\mathbf x-\mathbf y\rangle\\
&\le
-\left(\min_{1\le i\le m}d_i\right)
\|\mathbf x-\mathbf y\|^2
+\varepsilon L_G\|\mathbf x-\mathbf y\|^2\\
&=-\kappa\|\mathbf x-\mathbf y\|^2.
\end{aligned}
\]
By the differential equation for $\Phi$, we have
\[
\frac{d}{dt}
\|\Phi(t,\mathbf x)-\Phi(t,\mathbf y)\|^2
\le
-2\kappa\|\Phi(t,\mathbf x)-\Phi(t,\mathbf y)\|^2.
\]
By the product rule, we have
\[
\frac{d}{dt}
\left(
e^{2\kappa t}
\|\Phi(t,\mathbf x)-\Phi(t,\mathbf y)\|^2
\right)
\le0.
\]
Integrating from $0$ to $t$, we obtain
\begin{equation}\label{eq:pair-contraction}
\|\Phi(t,\mathbf x)-\Phi(t,\mathbf y)\|
\le e^{-\kappa t}\|\mathbf x-\mathbf y\|
\qquad
(t\ge0,\ \mathbf x,\mathbf y\in\mathbb R^m).
\end{equation}

The vector field $b_0(\mathbf x)=-D\mathbf x$ generates the flow
\[
T(t,\mathbf x)
=
(e^{-d_1t}x_1,\ldots,e^{-d_mt}x_m).
\]
By the first hypothesis of the theorem, we have $T(t,\mathbf x)\in\operatorname{int}C$ for all $t>0$ and $\mathbf x\in C$.
Since $\delta_C$ vanishes on $\mathbb R^m\setminus\operatorname{int}C$, we have $b=b_0$ on this set.
By Lemma~\ref{lem:boundary-robustness}, we have
\begin{equation}\label{eq:matrix-C-interior}
\Phi(t,\mathbf x)\in\operatorname{int}C
\qquad(t>0,\ \mathbf x\in C).
\end{equation}

Fix $\mathbf x\in X$ and $t>0$, and define
\[
\lambda(s)=\varepsilon\delta_C(\Phi(s,\mathbf x))
\qquad(0\le s\le t).
\]
The function $\lambda$ is continuous, and
\[
\frac{d}{ds}\Phi(s,\mathbf x)
=
-D\Phi(s,\mathbf x)
+\lambda(s)\bigl(B\Phi(s,\mathbf x)+\mathbf c\bigr).
\]
By~\eqref{eq:matrix-C-interior}, we have $\lambda(s)>0$ for $0<s\le t$.
Since $\lambda(0)$ may be zero, we first prove that $\Phi(t/2,\mathbf x)\in Y$.

For $0<r<t/2$, consider the initial value problem
\[
\mathbf u_r'(s)
=
-D\mathbf u_r(s)
+\lambda(s)\bigl(B\mathbf u_r(s)+\mathbf c\bigr),
\qquad
\mathbf u_r(r)=\mathbf x,
\]
with the same coefficient function $\lambda$.
By~\cite[Theorem~2.2 and Corollary~2.6]{Teschl}, there exists a unique solution $\mathbf u_r:[r,t/2]\to\mathbb R^m$.
The function $\lambda$ is positive on $[r,t/2]$, and $\mathbf u_r(r)=\mathbf x\in Y$.
By the second hypothesis, after translating the initial time to zero, we have $\mathbf u_r(t/2)\in\operatorname{int}Y$.
Set
\[
K_t=\max_{0\le s\le t}
\|-D+\lambda(s)B\|_{\mathrm{op}}.
\]
The number $K_t$ is finite because $\lambda$ is continuous.
By the equations for $\mathbf u_r$ and $\Phi$, we have
\[
\frac{d}{ds}
\bigl(\mathbf u_r(s)-\Phi(s,\mathbf x)\bigr)
=
\bigl(-D+\lambda(s)B\bigr)
\bigl(\mathbf u_r(s)-\Phi(s,\mathbf x)\bigr).
\]
Consequently,
\[
\frac12\frac{d}{ds}
\|\mathbf u_r(s)-\Phi(s,\mathbf x)\|^2
\le K_t\|\mathbf u_r(s)-\Phi(s,\mathbf x)\|^2.
\]
By the product rule, we have
\[
\frac{d}{ds}
\left(
e^{-2K_t s}
\|\mathbf u_r(s)-\Phi(s,\mathbf x)\|^2
\right)
\le0.
\]
Integrating from $r$ to $t/2$, we obtain
\[
\|\mathbf u_r(t/2)-\Phi(t/2,\mathbf x)\|
\le e^{K_t(t/2-r)}
\|\mathbf x-\Phi(r,\mathbf x)\|.
\]
By continuity of $\Phi$, the right-hand side tends to zero as $r\to0^+$.
Since $\mathbf u_r(t/2)\in Y$ and $Y$ is closed, we have $\Phi(t/2,\mathbf x)\in Y$.
By the second hypothesis on $[t/2,t]$, where $\lambda$ is positive, we have $\Phi(t,\mathbf x)\in\operatorname{int}Y$.
Together with~\eqref{eq:matrix-C-interior}, we obtain
\begin{equation}\label{eq:inward}
\Phi(t,\mathbf x)
\in\operatorname{int}C\cap\operatorname{int}Y
=\operatorname{int}X
\qquad(t>0,\ \mathbf x\in X).
\end{equation}

Choose $\mathbf x_0\in X$ and $R>0$ such that $X\subseteq R\overline B^m$.
By~\eqref{eq:inward}, we have $\Phi(s,\mathbf x_0)\in X$ for all $s\ge0$.
By the flow identity and~\eqref{eq:pair-contraction}, we have
\[
\begin{aligned}
\|\Phi(t+s,\mathbf x_0)-\Phi(t,\mathbf x_0)\|
&=
\|\Phi(t,\Phi(s,\mathbf x_0))-\Phi(t,\mathbf x_0)\|\\
&\le e^{-\kappa t}
\|\Phi(s,\mathbf x_0)-\mathbf x_0\|\\
&\le 2R e^{-\kappa t}
\qquad(s,t\ge0).
\end{aligned}
\]
By the Cauchy criterion, the limit $\mathbf p=\lim_{t\to\infty}\Phi(t,\mathbf x_0)$ exists in $\mathbb R^m$.
Since $X$ is closed, we have $\mathbf p\in X$.
For every $s\in\mathbb R$, by continuity and the flow identity, we have
\[
\Phi(s,\mathbf p)
=
\lim_{t\to\infty}\Phi(s,\Phi(t,\mathbf x_0))
=
\lim_{t\to\infty}\Phi(s+t,\mathbf x_0)
=
\mathbf p.
\]
In particular, by~\eqref{eq:inward}, we have $\mathbf p=\Phi(1,\mathbf p)\in\operatorname{int}X$.

For every $\mathbf x\in X$, by~\eqref{eq:inward}, we have $\Phi(t,\mathbf x)\in\operatorname{int}X$ for $t>0$.
By continuity of $\Phi$, we also have $\Phi(t,\mathbf x)\to\mathbf x$ as $t\to0^+$.
Therefore $X\subseteq\overline{\operatorname{int}X}$.
The reverse inclusion holds because $X$ is closed, so
\[
X=\overline{\operatorname{int}X}.
\]
Define
\[
U=\{\mathbf x-\mathbf p:\mathbf x\in\operatorname{int}X\}
\]
and $F:\mathbb R\times\mathbb R^m\to\mathbb R^m$ by
\[
F(t,\mathbf z)=\Phi(t,\mathbf p+\mathbf z)-\mathbf p.
\]
Then $U$ is nonempty, bounded, and open, and $\overline U=\{\mathbf x-\mathbf p:\mathbf x\in X\}$.
The map $F$ is continuous and satisfies the flow identities.
Since $\Phi(t,\mathbf p)=\mathbf p$, by~\eqref{eq:pair-contraction}, we have
\[
\|F(t,\mathbf z)\|
=
\|\Phi(t,\mathbf p+\mathbf z)-\Phi(t,\mathbf p)\|
\le e^{-\kappa t}\|\mathbf z\|
\qquad(t\ge0,\ \mathbf z\in\mathbb R^m).
\]
For $t>0$ and $\mathbf z\ne\mathbf0$, we have $e^{-\kappa t}\|\mathbf z\|<\|\mathbf z\|$.
Hence $F$ is contractive.
By~\eqref{eq:inward}, we also have
\[
F(t,\mathbf z)\in U
\qquad(t>0,\ \mathbf z\in\overline U).
\]
By Lemma~\ref{lem:gkl}, there exists a homeomorphism $g:\overline U\to\overline B^m$ such that $g(\partial U)=S^{m-1}$.
Since $U$ is open and $X$ is closed, we have
\[
\partial U
=
\overline U\setminus U
=
\{\mathbf x-\mathbf p:\mathbf x\in\partial X\}.
\]
Define $h:X\to\overline B^m$ by $h(\mathbf x)=g(\mathbf x-\mathbf p)$.
Then $h$ is a homeomorphism and $h(\partial X)=S^{m-1}$.
\end{proof}

\section{Lorentzian and stable polynomials on cones}
\label{sec:cones}

In this section, we prove Theorems~\ref{thm:lor-main} and~\ref{thm:stable-main}.

\subsection{Preliminaries and reduction by lineality}
\label{subsec:cone-preliminaries}

For $d\ge0$ and $f\in\mathrm H_d(\mathbb R^n)$, define
\[
\|f\|=\|\operatorname{coef}_d(f)\|,
\]
where the norm on the right-hand side is the Euclidean norm on $\mathbb R^M$.
For $f,g\in\mathrm H_d(\mathbb R^n)$, we have $\rho(f,g)=\|f-g\|$.
An \emph{affine subspace} of $\mathrm H_d(\mathbb R^n)$ is a set
\[
\mathcal A=f_0+E
=
\{f_0+g:g\in E\},
\]
where $f_0\in\mathrm H_d(\mathbb R^n)$ and $E\subseteq\mathrm H_d(\mathbb R^n)$ is a linear subspace.
The space $E$ is called the \emph{translation space} of $\mathcal A$, and $\dim\mathcal A=\dim E$.
We equip $\mathcal A$ with the restriction of the metric $\rho$.
For $X\subseteq\mathcal A$, let $\operatorname{int}_{\mathcal A}X$ and $\overline X^{\mathcal A}$ denote the interior and closure of $X$ in $\mathcal A$, respectively.
The boundary of $X$ in $\mathcal A$ is $\partial_{\mathcal A}X=\overline X^{\mathcal A}\setminus\operatorname{int}_{\mathcal A}X$.

For $d\ge1$ and $f\in\mathrm H_d(\mathbb R^n)$, the \emph{symmetric polarization} of $f$ is the map $P_f:(\mathbb R^n)^d\to\mathbb R$ defined by
\begin{equation}\label{eq:polarization}
P_f(\mathbf u_1,\ldots,\mathbf u_d)
=
\frac{1}{d!}
D_{\mathbf u_1}\cdots D_{\mathbf u_d}f.
\end{equation}
The expression on the right-hand side is constant because $f$ is homogeneous of degree $d$.
For $\alpha,\beta\in\mathbb R$ and $\mathbf a,\mathbf b\in\mathbb R^n$, we have $D_{\alpha\mathbf a+\beta\mathbf b}=\alpha D_{\mathbf a}+\beta D_{\mathbf b}$.
Since $D_{\mathbf a}D_{\mathbf b}=D_{\mathbf b}D_{\mathbf a}$, $P_f$ is a symmetric $d$-linear form.
By homogeneity, we have
\[
P_f(\mathbf x,\ldots,\mathbf x)
=
\frac{1}{d!}D_{\mathbf x}^{\,d}f
=
f(\mathbf x)
\qquad(\mathbf x\in\mathbb R^n).
\]

For a set $B$, let $B^0$ denote the singleton consisting of the empty tuple.
For $d\ge2$ and $\mathbf u=(\mathbf u_3,\ldots,\mathbf u_d)$, define the symmetric bilinear form
\begin{equation}\label{eq:mixed-form}
Q_f(\mathbf a,\mathbf b;\mathbf u)
=P_f(\mathbf a,\mathbf b,\mathbf u_3,\ldots,\mathbf u_d).
\end{equation}
Let $\mathbf e_1,\ldots,\mathbf e_n$ be the standard basis of $\mathbb R^n$.
By~\eqref{eq:polarization}, we have
\[
\operatorname{Hess}(D_{\mathbf u_3}\cdots D_{\mathbf u_d}f)
=d!\bigl(Q_f(\mathbf e_i,\mathbf e_j;\mathbf u)\bigr)_{i,j=1}^n.
\]
When $d=2$, the tuple $\mathbf u$ is empty and the matrix on the left is $\operatorname{Hess}(f)$.

For a symmetric bilinear form $Q$ on $\mathbb R^k$, let $n_+(Q)$ denote the number of positive eigenvalues of a matrix representing $Q$, counted with multiplicity.
The form $Q$ is called \emph{positive definite} if $Q(\mathbf a,\mathbf a)>0$ for every nonzero $\mathbf a\in\mathbb R^k$, and \emph{negative definite} if $Q(\mathbf a,\mathbf a)<0$ for every nonzero $\mathbf a\in\mathbb R^k$.
It is called \emph{nonpositive} if $Q(\mathbf a,\mathbf a)\le0$ for every $\mathbf a\in\mathbb R^k$.
We use the following standard linear algebra facts.

\begin{proposition}[{\cite[Chapter~4]{HornJohnson}}]
\label{prop:inertia-facts}
Let $Q$ be a symmetric bilinear form on $\mathbb R^k$.
\begin{enumerate}
\item
The number $n_+(Q)$ is independent of the basis used to represent $Q$.

\item
The eigenvalues of a real symmetric matrix, ordered from largest to smallest and counted with multiplicity, depend continuously on the matrix entries.

\item
If $k\ge2$, then $n_+(Q)\le1$ if and only if $Q$ is nonpositive on a linear subspace of dimension $k-1$.

\item
If $k\ge2$, then the second largest eigenvalue of a
matrix representing $Q$ is negative if and only if $Q$
is negative definite on a linear subspace of dimension $k-1$.
\end{enumerate}
\end{proposition}

A nonzero polynomial in one real variable is called \emph{real-rooted} if all complex roots are real; nonzero constants are included.
A polynomial is called \emph{monic} if the leading coefficient is one.
A root $a$ has \emph{multiplicity} $j$ if $p(z)=(z-a)^j q(z)$ with $q(a)\ne0$; a root of multiplicity one is called \emph{simple}.
We use the following closure property, which is the univariate real-coefficient form of the stability closure property in~\cite[Proposition~3.1(3)]{BBL}.

\begin{theorem}[{\cite[Proposition~3.1(3)]{BBL}}]
\label{thm:real-rooted-limits}
Let $(p_j)_{j\ge1}$ be a sequence of real-rooted real polynomials whose degrees have a common finite upper bound.
If the coefficients of $p_j$ converge to the coefficients of a nonzero polynomial $p$, then $p$ is real-rooted.
In particular, the set of monic real-rooted polynomials of a fixed degree is closed in the corresponding coefficient space.
\end{theorem}

A topological space is called \emph{connected} if it cannot be expressed as the union of two disjoint nonempty open subsets.
A \emph{connected component} is a maximal connected subset.
For $d\ge1$, a polynomial $f\in\mathrm H_d(\mathbb R^n)$ is called \emph{hyperbolic with respect to $\mathbf v$} if $f(\mathbf v)>0$ and the polynomial $t\mapsto f(\mathbf x+t\mathbf v)$ is real-rooted for every $\mathbf x\in\mathbb R^n$.
The \emph{hyperbolicity cone} $\Gamma_{\mathbf v}(f)$ is the connected component of $\{\mathbf x\in\mathbb R^n:f(\mathbf x)\ne0\}$ containing $\mathbf v$.
The polynomial $f$ is called \emph{strictly hyperbolic with respect to $\mathbf v$} if it is hyperbolic with respect to $\mathbf v$ and $t\mapsto f(\mathbf x+t\mathbf v)$ has $d$ distinct real roots whenever $\mathbf x\notin\mathbb R\mathbf v$.

\begin{theorem}[G\aa rding~\cite{Garding}; see also {\cite[Proposition~4.4]{BBL}}]
\label{thm:garding}
Suppose that $f$ is hyperbolic with respect to $\mathbf v$.
Then $\Gamma_{\mathbf v}(f)$ is an open convex cone.
For every $\mathbf u\in\Gamma_{\mathbf v}(f)$, the polynomial $f$ is hyperbolic with respect to $\mathbf u$, and $\Gamma_{\mathbf u}(f)=\Gamma_{\mathbf v}(f)$.
\end{theorem}

\begin{corollary}\label{cor:stable-hyperbolic}
Let $d\ge1$ and $f\in\mathrm H_d(\mathbb R^n)$ satisfy $f(\mathbf v)=1$.
Then $f$ is $K$-stable if and only if $f$ is hyperbolic with respect to $\mathbf v$ and $K\subseteq\Gamma_{\mathbf v}(f)$.
\end{corollary}

\begin{proof}
Suppose that $f$ is $K$-stable.
For $\mathbf y\in K$, by homogeneity, we have $f(\mathrm i\mathbf y)=\mathrm i^d f(\mathbf y)\ne0$.
Since $K$ is convex and contains $\mathbf v$, it is a connected subset of $\{f\ne0\}$ containing $\mathbf v$.
Therefore $K\subseteq\Gamma_{\mathbf v}(f)$.
Fix $\mathbf x\in\mathbb R^n$.
The polynomial $t\mapsto f(\mathbf x+t\mathbf v)$ has degree $d$ and leading coefficient $f(\mathbf v)=1$.
If $a+\mathrm i b$ were a root with $b>0$, then
\[
f((\mathbf x+a\mathbf v)+\mathrm i b\mathbf v)=0,
\]
contrary to $K$-stability, because $b\mathbf v\in K$.
Since the coefficients are real, every nonreal root occurs with the complex conjugate root.
Consequently all roots are real, and $f$ is hyperbolic with respect to $\mathbf v$.

Conversely, suppose that $f$ is hyperbolic with respect to $\mathbf v$ and $K\subseteq\Gamma_{\mathbf v}(f)$.
Fix $\mathbf x\in\mathbb R^n$ and $\mathbf y\in K$.
By Theorem~\ref{thm:garding}, the polynomial $t\mapsto f(\mathbf x+t\mathbf y)$ is real-rooted and has nonzero leading coefficient $f(\mathbf y)$.
It does not vanish at $t=\mathrm i$.
Hence $f(\mathbf x+\mathrm i\mathbf y)\ne0$.
\end{proof}

We also use the following known relation between hyperbolic and Lorentzian polynomials; see the discussion following Example~2.2 in~\cite[Section~2]{BL}.

\begin{proposition}[{\cite[Section~2]{BL}}]
\label{prop:hyperbolic-lorentzian}
A homogeneous polynomial hyperbolic with respect to $\mathbf v$ is $\Gamma_{\mathbf v}(f)$-Lorentzian.
\end{proposition}

\begin{corollary}\label{cor:stable-lorentzian}
For every $d\ge0$, we have $\mathcal S_{K,d,\mathbf v}\subseteq\mathcal L_{K,d,\mathbf v}$.
\end{corollary}

\begin{proof}
In degree zero both spaces consist of the constant polynomial $1$.
Let $d\ge1$ and $f\in\mathcal S_{K,d,\mathbf v}$.
By Corollary~\ref{cor:stable-hyperbolic}, the polynomial $f$ is hyperbolic with respect to $\mathbf v$ and $K\subseteq\Gamma_{\mathbf v}(f)$.
By Proposition~\ref{prop:hyperbolic-lorentzian}, the defining conditions for a Lorentzian polynomial hold for directions in $\Gamma_{\mathbf v}(f)$, and hence for directions in $K$.
Since $f(\mathbf v)=1$, we have $f\in\mathcal L_{K,d,\mathbf v}$.
\end{proof}

Recall that $S=\overline K\cap(-\overline K)$ and $r=n-\dim S$.
For $d\ge1$, define the affine subspace
\begin{equation}\label{eq:original-affine}
\mathcal A_{K,d,\mathbf v}
=\{f\in\mathrm H_d(\mathbb R^n):
D_{\mathbf w}f=0\text{ for every }\mathbf w\in S,\ f(\mathbf v)=1\}.
\end{equation}

\begin{proposition}[{\cite[Propositions~2.6 and~2.7]{BL}}]
\label{prop:BL-input}
Let $d\ge1$.
If $f$ is $K$-Lorentzian and $A:\mathbb R^p\to\mathbb R^n$ is linear with $A^{-1}(K)\ne\varnothing$, then $f\circ A$ is $A^{-1}(K)$-Lorentzian.
If $f\in\mathrm H_d(\mathbb R^n)$ satisfies $P_f>0$ on $K^d$, then $D_{\mathbf w}f=0$ for every $\mathbf w\in S$.
\end{proposition}

We use finite-dimensional separation of convex sets in the following form: if $O\subseteq\mathbb R^n$ is nonempty, open, and convex and $\mathbf x\notin O$, then there exists a nonzero linear functional $a$ such that $a(\mathbf u)\le a(\mathbf x)$ for all $\mathbf u\in O$; see~\cite[Section~11]{Rockafellar}.

\begin{corollary}\label{cor:quotient}
There exist a nonempty open convex cone $\widehat K\subsetneq\mathbb R^r$ with
$\overline{\widehat K}\cap(-\overline{\widehat K})=\{\mathbf0\}$, a point $\widehat{\mathbf v}\in\widehat K$, and a surjective linear map $\pi:\mathbb R^n\to\mathbb R^r$ such that pullback by $\pi$ is an affine homeomorphism
\[
\{\widehat f\in\mathrm H_d(\mathbb R^r):\widehat f(\widehat{\mathbf v})=1\}
\longrightarrow\mathcal A_{K,d,\mathbf v}.
\]
This homeomorphism restricts to homeomorphisms
\[
\mathcal L_{\widehat K,d,\widehat{\mathbf v}}
\longrightarrow\mathcal L_{K,d,\mathbf v},
\qquad
\mathcal S_{\widehat K,d,\widehat{\mathbf v}}
\longrightarrow\mathcal S_{K,d,\mathbf v}.
\]
In particular, both spaces on $K$ are contained in $\mathcal A_{K,d,\mathbf v}$, and
\begin{equation}\label{eq:effective-dimension}
\dim\mathcal A_{K,d,\mathbf v}=\binom{r+d-1}{d}-1.
\end{equation}
\end{corollary}

\begin{proof}
Set $P=\overline K$.
The set $P$ is a closed convex cone.
Since $S=P\cap(-P)$ is closed under addition and multiplication by real scalars, it is a linear subspace.
We also have $P+S=P$.

We first verify that $K=\operatorname{int}P$.
The inclusion $K\subseteq\operatorname{int}P$ follows because $K$ is open.
Suppose that $\mathbf x\in\operatorname{int}P\setminus K$.
By separation, there is a nonzero linear functional $a$ such that $a(\mathbf u)\le a(\mathbf x)$ for all $\mathbf u\in K$.
By continuity, the same inequality holds for all $\mathbf u\in P$.
Choose $\mathbf z$ with $a(\mathbf z)>0$.
For sufficiently small $t>0$, we have $\mathbf x+t\mathbf z\in P$ and $a(\mathbf x+t\mathbf z)>a(\mathbf x)$, a contradiction.
Therefore $K=\operatorname{int}P$, and $P+S=P$ yields $K+S=K$.

Let $S^\perp=\{\mathbf x\in\mathbb R^n:\langle\mathbf x,\mathbf s\rangle=0\text{ for every }\mathbf s\in S\}$ be the orthogonal complement of $S$.
Choose an orthonormal basis of $S^\perp$, and let $\pi:\mathbb R^n\to\mathbb R^r$ be the orthogonal projection written in this basis.
Under this identification, $\pi(P)=P\cap S^\perp$ and $\pi(K)=K\cap S^\perp$.
Indeed, subtracting the component in $S$ leaves membership in $P$ and in $K$ unchanged.
Set $\widehat K=\pi(K)$ and $\widehat{\mathbf v}=\pi(\mathbf v)$.
Then $\widehat K$ is open and convex, and $\pi^{-1}(\widehat K)=K$.
The set $\pi(P)$ is closed.
For every $\mathbf p\in P$, a sequence of points of $K$ converging to $\mathbf p$ projects to a sequence in $\widehat K$ converging to $\pi(\mathbf p)$.
Consequently $\overline{\widehat K}=\pi(P)$.
If both $\mathbf u$ and $-\mathbf u$ belong to $\pi(P)=P\cap S^\perp$, then $\mathbf u\in S\cap S^\perp=\{\mathbf0\}$.
Thus $\overline{\widehat K}$ contains no nonzero linear subspace.
Since $K\ne\mathbb R^n$, we have $r\ge1$.
Since $r\ge1$ and $\overline{\widehat K}\cap(-\overline{\widehat K})=\{\mathbf0\}$, we have $\widehat K\ne\mathbb R^r$.
It does not contain $\mathbf0$: an open cone containing $\mathbf0$ contains a ball centered at $\mathbf0$ and therefore, by positive scaling, equals the whole ambient space.
In particular, $\widehat{\mathbf v}\ne\mathbf0$.

By Proposition~\ref{prop:BL-input}, every $f\in\mathcal L_{K,d,\mathbf v}$ satisfies $D_{\mathbf w}f=0$ for $\mathbf w\in S$.
By Corollary~\ref{cor:stable-lorentzian}, the same assertion holds for $f\in\mathcal S_{K,d,\mathbf v}$.
More generally, $D_{\mathbf w}f=0$ for every $\mathbf w\in S$ holds if and only if $f$ is constant on every affine subspace $\mathbf x+S$.
In coordinates adapted to $S\oplus S^\perp$, this is equivalent to independence of the coordinates in $S$.
Such a polynomial has a unique representation $f=\widehat f\circ\pi$ with $\widehat f\in\mathrm H_d(\mathbb R^r)$.

By Proposition~\ref{prop:BL-input}, pullback by $\pi$ preserves the Lorentzian property.
Conversely, let $\iota:\mathbb R^r\to\mathbb R^n$ be the inclusion with image $S^\perp$ in the chosen coordinates.
Then $\iota^{-1}(K)=\widehat K$ and $\widehat f=f\circ\iota$, so the same proposition proves the converse.
For stability, extend $\pi$ complex linearly.
The identity $f(\mathbf x+\mathrm i\mathbf y)=\widehat f(\pi\mathbf x+\mathrm i\pi\mathbf y)$ and the equality
$\pi(\mathbb R^n+\mathrm iK)=\mathbb R^r+\mathrm i\widehat K$
prove equivalence of the two stability conditions.
Normalization is preserved because $f(\mathbf v)=\widehat f(\widehat{\mathbf v})$.

Pullback is a linear isomorphism from $\mathrm H_d(\mathbb R^r)$ onto the subspace of polynomials annihilated by all $D_{\mathbf w}$ with $\mathbf w\in S$.
A linear isomorphism between finite-dimensional real vector spaces is a homeomorphism in the coefficient topologies.
Evaluation at $\widehat{\mathbf v}\ne\mathbf0$ is a nonzero linear functional on $\mathrm H_d(\mathbb R^r)$.
The level set at $1$ therefore has dimension $\binom{r+d-1}{d}-1$.
\end{proof}

In degree zero, both normalized spaces consist of the constant polynomial $1$, so both asserted dimensions are zero.
For the remainder of this section, let $d\ge1$.
By Corollary~\ref{cor:quotient}, it suffices to work in $\mathbb R^r$ with the reduced cone and point.
We write $K$ and $\mathbf v$ for $\widehat K$ and $\widehat{\mathbf v}$.
Thus $P=\overline K$ is \emph{pointed}, meaning that $P\cap(-P)=\{\mathbf0\}$.

\subsection{A common compact set and diagonal contraction}
\label{subsec:common-contraction}

We use two standard finite-dimensional facts: the convex hull of a compact set is compact, and a compact convex set not containing $\mathbf0$ can be strictly separated from $\mathbf0$ by a linear functional; see~\cite[Sections~11 and~17]{Rockafellar}.

The \emph{dual cone} of $P$ is
\[
P^*=\{\ell:\mathbb R^r\to\mathbb R\text{ linear}:\ell(\mathbf x)\ge0\text{ for every }\mathbf x\in P\}.
\]
Interior in the space of linear functionals is taken in the Euclidean topology of the coefficient vectors.
Every $\ell\in\operatorname{int}P^*$ is strictly positive on $P\setminus\{\mathbf0\}$.
Indeed, if $\ell(\mathbf u)=0$ for a nonzero $\mathbf u\in P$, choose a linear functional $a$ with $a(\mathbf u)>0$.
Then $\ell-\varepsilon a\notin P^*$ for every $\varepsilon>0$, contrary to $\ell\in\operatorname{int}P^*$.

\begin{lemma}\label{lem:cone-base}
There exists $\ell\in\operatorname{int}P^*$ with $\ell(\mathbf v)=1$.
For such a functional, the set
\begin{equation}\label{eq:cone-base}
\Delta=\{\mathbf u\in P:\ell(\mathbf u)=1\}
\end{equation}
is compact, spans $\mathbb R^r$, and contains $\mathbf v$.
There exists $\delta>0$ such that $\mathbf v-\delta\mathbf u\in P$ for every $\mathbf u\in\Delta$.
\end{lemma}

\begin{proof}
The set $P\cap S^{r-1}$ is nonempty and compact.
The convex hull of this set does not contain $\mathbf0$.
Otherwise, after removing zero coefficients, we would have
$\sum_{j=1}^k a_j\mathbf u_j=\mathbf0$ with $a_j>0$, $\sum_j a_j=1$, and $\mathbf u_j\in P\cap S^{r-1}$.
The case $k=1$ is impossible.
For $k\ge2$, we would have $-\mathbf u_1=\sum_{j=2}^k(a_j/a_1)\mathbf u_j\in P$, contrary to pointedness.
By compactness and strict separation, there exist a linear functional $\ell_0$ and $a>0$ such that
\[
\ell_0(\mathbf u)\ge a\qquad(\mathbf u\in P\cap S^{r-1}).
\]
By positive scaling, $\ell_0(\mathbf u)\ge a\|\mathbf u\|$ for every $\mathbf u\in P$.
If the Euclidean norm of the coefficients of a linear functional $b$ is smaller than $a$, then $|b(\mathbf u)|<a\|\mathbf u\|$ for $\mathbf u\ne\mathbf0$.
Consequently $\ell_0+b\in P^*$, so $\ell_0\in\operatorname{int}P^*$.
Since $\mathbf v\in P\setminus\{\mathbf0\}$, we may set $\ell=\ell_0/\ell_0(\mathbf v)$.

The same lower bound shows that $\Delta$ is bounded.
It is closed, so it is compact.
Every nonzero vector of $P$ is a positive multiple of a point of $\Delta$.
Since $P$ contains the nonempty open set $K$, it spans $\mathbb R^r$; hence $\Delta$ also spans $\mathbb R^r$.
We have $\mathbf v\in\Delta$ by normalization.
Choose $\eta>0$ such that $\{\mathbf x:\|\mathbf x-\mathbf v\|<\eta\}\subseteq K$.
Since $\Delta$ is bounded, there is $\delta>0$ with $\delta\|\mathbf u\|<\eta$ for every $\mathbf u\in\Delta$.
Then $\mathbf v-\delta\mathbf u\in K\subseteq P$.
\end{proof}

Fix $\ell$ and $\Delta$ as in Lemma~\ref{lem:cone-base}, and set
\begin{equation}\label{eq:reduced-affine}
\mathcal A=\{f\in\mathrm H_d(\mathbb R^r):f(\mathbf v)=1\},
\qquad
E=\{g\in\mathrm H_d(\mathbb R^r):g(\mathbf v)=0\}.
\end{equation}
The translation space of $\mathcal A$ is $E$.
Define
\begin{equation}\label{eq:common-compact}
\mathcal C=\{f\in\mathcal A:P_f(\mathbf u_1,\ldots,\mathbf u_d)\ge0
\text{ for all }\mathbf u_1,\ldots,\mathbf u_d\in\Delta\}.
\end{equation}

\begin{lemma}\label{lem:common-compact}
The set $\mathcal C$ is nonempty, compact, and convex.
For every $f\in\mathcal C$, we have $P_f>0$ on $K^d$.
Moreover, $\mathcal L_{K,d,\mathbf v}\subseteq\mathcal C$ and $\mathcal S_{K,d,\mathbf v}\subseteq\mathcal C$.
\end{lemma}

\begin{proof}
Each condition in~\eqref{eq:common-compact} is a closed affine inequality in the coefficients of $f$.
Therefore $\mathcal C$ is closed and convex in $\mathcal A$.
The polynomial $\ell^d$ belongs to $\mathcal C$, so $\mathcal C$ is nonempty.
For $f\in\mathcal C$, multilinearity and positive scaling yield $P_f\ge0$ on $P^d$.
This includes tuples with a zero vector, because the corresponding value is zero.

Choose $\delta$ as in Lemma~\ref{lem:cone-base}.
For $\mathbf u_1,\ldots,\mathbf u_d\in\Delta$, substitute $\mathbf v=\delta\mathbf u_i+(\mathbf v-\delta\mathbf u_i)$ in the $i$th argument of $P_f(\mathbf v,\ldots,\mathbf v)$.
Every term in the resulting multilinear expansion is nonnegative.
Consequently
\[
0\le P_f(\mathbf u_1,\ldots,\mathbf u_d)\le\delta^{-d}f(\mathbf v)=\delta^{-d}.
\]
Choose a basis $\mathbf b_1,\ldots,\mathbf b_r$ from $\Delta$.
For nonnegative integers $\alpha_1,\ldots,\alpha_r$ with $\sum_j\alpha_j=d$, the coefficient of $t_1^{\alpha_1}\cdots t_r^{\alpha_r}$ in $f(\sum_j t_j\mathbf b_j)$ is
\[
\frac{d!}{\alpha_1!\cdots\alpha_r!}
P_f(\underbrace{\mathbf b_1,\ldots,\mathbf b_1}_{\alpha_1},\ldots,
\underbrace{\mathbf b_r,\ldots,\mathbf b_r}_{\alpha_r}).
\]
The preceding estimate bounds all these coefficients uniformly over $f\in\mathcal C$.
The change from this basis to the standard coordinates is a fixed linear isomorphism on the polynomial space.
Hence $\mathcal C$ is bounded in the coefficient norm, and it is compact.

Fix $\mathbf u_1,\ldots,\mathbf u_d\in K$.
For each $i$, choose $a_i>0$ sufficiently small that $\mathbf u_i-a_i\mathbf v\in K$.
By multilinearity and nonnegativity on $P^d$, we have
\[
P_f(\mathbf u_1,\ldots,\mathbf u_d)
\ge\left(\prod_{i=1}^d a_i\right)P_f(\mathbf v,\ldots,\mathbf v)
=\prod_{i=1}^d a_i>0.
\]
If $f\in\mathcal L_{K,d,\mathbf v}$, continuity of $P_f$ and $P=\overline K$ yield $P_f\ge0$ on $P^d$, so $f\in\mathcal C$.
The assertion for stable polynomials follows by Corollary~\ref{cor:stable-lorentzian}.
\end{proof}

For $f\in\mathcal A$, define
\[
\mu(f)=\min_{(\mathbf u_1,\ldots,\mathbf u_d)\in\Delta^d}
P_f(\mathbf u_1,\ldots,\mathbf u_d).
\]
This minimum exists by compactness.
The map $(g,\mathbf u_1,\ldots,\mathbf u_d)\mapsto P_g(\mathbf u_1,\ldots,\mathbf u_d)$ is continuous and linear in $g$.
The unit sphere of the coefficient norm and the set $\Delta^d$ are compact.
There is therefore a finite constant $M_0>0$ such that
\[
|P_g(\mathbf u_1,\ldots,\mathbf u_d)|\le M_0\|g\|
\qquad(g\in\mathrm H_d(\mathbb R^r),\ (\mathbf u_1,\ldots,\mathbf u_d)\in\Delta^d).
\]
Applying this inequality to $g=f_1-f_2$ and interchanging $f_1,f_2$, we obtain
\begin{equation}\label{eq:mu-lipschitz}
|\mu(f_1)-\mu(f_2)|\le M_0\|f_1-f_2\|
\qquad(f_1,f_2\in\mathcal A).
\end{equation}
In particular, $\mu(f)>0$ implies $f\in\operatorname{int}_{\mathcal A}\mathcal C$.

Set $f_*=\ell^d$.
For $\mathbf y\in K$, the complex number $\ell(\mathbf x+\mathrm i\mathbf y)$ has positive imaginary part.
Therefore $f_*$ is $K$-stable, and $f_*(\mathbf v)=1$.
By Corollary~\ref{cor:stable-lorentzian}, $f_*$ also belongs to $\mathcal L_{K,d,\mathbf v}$.
If $r=1$, the affine space $\mathcal A$ consists only of $f_*$.
Both normalized spaces are then one-point spaces, as asserted in Theorems~\ref{thm:lor-main} and~\ref{thm:stable-main}.
We assume $r\ge2$ below.

Let $W=\ker\ell$ and choose an orthonormal basis of $W$ for the restricted Euclidean inner product.
Let $\mathrm H_k(W)$ denote the polynomials homogeneous of degree $k$ in these coordinates.
Every $\mathbf x\in\mathbb R^r$ has a unique representation $\mathbf x=\mathbf y+s\mathbf v$ with $\mathbf y\in W$ and $s=\ell(\mathbf x)$.
Every $f\in\mathcal A$ has a unique expansion
\begin{equation}\label{eq:graded-coordinates}
f(\mathbf y+s\mathbf v)=s^d+\sum_{k=1}^d s^{d-k}p_k(\mathbf y),
\qquad p_k\in\mathrm H_k(W).
\end{equation}
For $g\in E$, the same expansion has zero coefficient of $s^d$.
Fix an ordering of the monomials in each $\mathrm H_k(W)$.
The ordered coefficient vectors of $p_1,\ldots,p_d$ define a linear isomorphism $J:E\to\mathbb R^N$, where
\[
N=\dim E=\binom{r+d-1}{d}-1.
\]
Define the affine homeomorphism $\Psi:\mathcal A\to\mathbb R^N$ by $\Psi(f)=J(f-f_*)$.
For $q\ge1$, let $I_q$ denote the $q\times q$ identity matrix.
Set $q_k=\dim\mathrm H_k(W)=\binom{r+k-2}{k}$ and
\begin{equation}\label{eq:cone-diagonal}
D=\operatorname{diag}(I_{q_1},2I_{q_2},\ldots,dI_{q_d}).
\end{equation}

Define $G:\mathbb R\times\mathbb R^r\to\mathbb R^r$ and $T:\mathbb R\times\mathrm H_d(\mathbb R^r)\to\mathrm H_d(\mathbb R^r)$ by
\begin{equation}\label{eq:base-flow}
G(t,\mathbf x)=e^{-t}\mathbf x+(1-e^{-t})\ell(\mathbf x)\mathbf v,
\qquad
T(t,f)(\mathbf x)=f(G(t,\mathbf x)).
\end{equation}
We have $G(t,\mathbf y+s\mathbf v)=e^{-t}\mathbf y+s\mathbf v$.
Consequently $G(s,G(t,\mathbf x))=G(s+t,\mathbf x)$, $G(t,\mathbf v)=\mathbf v$, and $\ell(G(t,\mathbf x))=\ell(\mathbf x)$.
The map $T$ is a flow, each map $f\mapsto T(t,f)$ is linear, and $T(t,\mathcal A)=\mathcal A$.
For a subset $Z$ of the polynomial space, we write $T(t,Z)=\{T(t,f):f\in Z\}$.
By~\eqref{eq:graded-coordinates}, $T(t,f)$ has coefficient polynomials $e^{-kt}p_k$.
Hence
\begin{equation}\label{eq:diagonal-coordinates}
\Psi(T(t,f))=e^{-tD}\Psi(f)
\qquad(t\in\mathbb R,\ f\in\mathcal A),
\end{equation}
where $e^{-tD}$ is the diagonal matrix with diagonal entries equal to the exponentials of the corresponding entries of $-tD$.

Define the linear operator $L:\mathrm H_d(\mathbb R^r)\to\mathrm H_d(\mathbb R^r)$ by
\begin{equation}\label{eq:cone-generator}
Lf=\ell D_{\mathbf v}f-df.
\end{equation}
By differentiating $f(t\mathbf x)=t^df(\mathbf x)$ at $t=1$, we have $\sum_i x_i\partial_i f(\mathbf x)=df(\mathbf x)$.
Differentiation of~\eqref{eq:base-flow} at $t=0$ therefore yields
$\left.\frac{d}{dt}T(t,f)\right|_{t=0}=Lf$.
By the flow identity, $\frac{d}{dt}T(t,f)=L(T(t,f))=T(t,Lf)$ for every real $t$.
Moreover, $Lf_*=0$ and $(Lf)(\mathbf v)=0$, so $L$ takes values in $E$.

\begin{lemma}\label{lem:common-contraction}
For $t>0$, we have $T(t,\mathcal C)\subseteq\operatorname{int}_{\mathcal A}\mathcal C$.
The set $C=\Psi(\mathcal C)$ is nonempty, closed, and bounded, and
\[
e^{-tD}C\subseteq\operatorname{int}C\qquad(t>0).
\]
\end{lemma}

\begin{proof}
For $t\ge0$ and $\mathbf u\in\Delta$, we have
$G(t,\mathbf u)=e^{-t}\mathbf u+(1-e^{-t})\mathbf v\in\Delta$.
By multilinearity, for $\mathbf u_1,\ldots,\mathbf u_d\in\Delta$, we have
\[
P_{T(t,f)}(\mathbf u_1,\ldots,\mathbf u_d)
=P_f(G(t,\mathbf u_1),\ldots,G(t,\mathbf u_d)).
\]
If $f\in\mathcal C$, every term obtained by substituting the formula for $G(t,\mathbf u_i)$ is nonnegative.
The term containing $\mathbf v$ in all $d$ arguments is $(1-e^{-t})^df(\mathbf v)=(1-e^{-t})^d$.
Thus $\mu(T(t,f))\ge(1-e^{-t})^d>0$ for $t>0$.
By~\eqref{eq:mu-lipschitz}, we have $T(t,f)\in\operatorname{int}_{\mathcal A}\mathcal C$.
The remaining assertions follow by Lemma~\ref{lem:common-compact} and~\eqref{eq:diagonal-coordinates}, because $\Psi$ is an affine homeomorphism.
\end{proof}

If $d=1$, both normalized spaces equal $\mathcal C$.
Indeed, membership in $\mathcal C$ is equivalent to nonnegativity of the linear functional on $P$ together with normalization at $\mathbf v$.
By Lemma~\ref{lem:common-compact}, this is equivalent to positivity on $K$ and normalization.
A real linear polynomial positive on $K$ is $K$-stable because the imaginary part of $f(\mathbf x+\mathrm i\mathbf y)$ is $f(\mathbf y)>0$.
Conversely, Corollary~\ref{cor:stable-lorentzian} proves the required positivity for a normalized $K$-stable linear polynomial.
Apply Theorem~\ref{thm:main-criterion} with $C=\Psi(\mathcal C)$, $Y=\mathbb R^N$, the matrix $D$ in~\eqref{eq:cone-diagonal}, $B=0$, and $\mathbf c=\mathbf0$.
The first hypothesis holds by Lemma~\ref{lem:common-contraction}, and the second holds because $Y$ is the entire ambient space.
This proves both ball conclusions in degree one.
Together with the cases $d=0$ and $r=1$, we may assume $r,d\ge2$ in the remaining subsections.

\subsection{Lorentzian polynomials}
\label{subsec:lorentzian-cones}

Define
\begin{equation}\label{eq:lorentzian-structural-set}
\mathcal Y_L=\{f\in\mathcal A:
n_+(Q_f(\cdot,\cdot;\mathbf u))\le1
\text{ for every }\mathbf u\in\Delta^{d-2}\}.
\end{equation}
For $f\in\mathcal A$ and $\mathbf u\in\Delta^{d-2}$, let $\theta_2(f,\mathbf u)$ denote the second largest eigenvalue of the matrix of $Q_f(\cdot,\cdot;\mathbf u)$ in the standard basis of $\mathbb R^r$.
By Proposition~\ref{prop:inertia-facts}, this is a continuous function of $f$ and $\mathbf u$.
Therefore $\mathcal Y_L$ is closed in $\mathcal A$.
If $\theta_2(f,\mathbf u)<0$ for every $\mathbf u\in\Delta^{d-2}$, then
\[
\max_{\mathbf u\in\Delta^{d-2}}\theta_2(f,\mathbf u)<0.
\]
By continuity and compactness of the parameter set, the same strict inequality holds for all polynomials in a sufficiently small neighborhood of $f$ in $\mathcal A$.
Consequently
\begin{equation}\label{eq:lorentzian-interior-test}
\theta_2(f,\mathbf u)<0\text{ for every }\mathbf u\in\Delta^{d-2}
\quad\Longrightarrow\quad f\in\operatorname{int}_{\mathcal A}\mathcal Y_L.
\end{equation}

\begin{lemma}\label{lem:lorentzian-intersection}
We have $\mathcal C\cap\mathcal Y_L=\mathcal L_{K,d,\mathbf v}$.
Moreover, $T(t,\mathcal Y_L)\subseteq\mathcal Y_L$ for every $t\ge0$.
\end{lemma}

\begin{proof}
Let $f\in\mathcal L_{K,d,\mathbf v}$.
By Lemma~\ref{lem:common-compact}, $f\in\mathcal C$.
Approximate any tuple of vectors in $P$ by tuples in $K$.
The corresponding mixed forms converge coefficientwise.
By Proposition~\ref{prop:inertia-facts}, the corresponding second largest eigenvalues converge and remain nonpositive in the limit.
Thus $f\in\mathcal Y_L$.

Conversely, let $f\in\mathcal C\cap\mathcal Y_L$.
By Lemma~\ref{lem:common-compact}, $P_f>0$ on $K^d$.
For $\mathbf u_3,\ldots,\mathbf u_d\in K$, divide each vector by the positive number $\ell(\mathbf u_i)$ to obtain a parameter in $\Delta$.
By multilinearity, the mixed form before division is a positive multiple of the mixed form after division.
It therefore has at most one positive eigenvalue.
Moreover,
$Q_f(\mathbf v,\mathbf v;\mathbf u_3,\ldots,\mathbf u_d)>0$,
so it has at least one positive eigenvalue.
Hence $f$ is $K$-Lorentzian.

For $t\ge0$, we have $G(t,\Delta)\subseteq\Delta$, and $G(t,\cdot)$ is invertible with inverse $G(-t,\cdot)$.
For all $\mathbf a,\mathbf b\in\mathbb R^r$, we have
\begin{equation}\label{eq:mixed-pullback}
Q_{T(t,f)}(\mathbf a,\mathbf b;\mathbf u)
=Q_f(G(t,\mathbf a),G(t,\mathbf b);
G(t,\mathbf u_3),\ldots,G(t,\mathbf u_d)).
\end{equation}
By Sylvester's law of inertia in Proposition~\ref{prop:inertia-facts}, the right-hand side has at most one positive eigenvalue when $f\in\mathcal Y_L$.
Thus $T(t,f)\in\mathcal Y_L$.
\end{proof}

Choose linearly independent $\ell_1,\ldots,\ell_r\in\operatorname{int}P^*$.
Such a choice is possible because $\operatorname{int}P^*$ is a nonempty open subset of the $r$-dimensional space of linear functionals.
Define
\begin{equation}\label{eq:positive-mixed-form-polynomial}
H=\frac{\sum_{j=1}^r\ell_j^d}{\sum_{j=1}^r\ell_j(\mathbf v)^d}.
\end{equation}
The denominator is positive, and $H\in\mathcal A$.
This construction is the positive-definite subtraction used in the proof of~\cite[Theorem~9.1]{BL}.
For every nonzero $\mathbf a\in\mathbb R^r$ and $\mathbf u\in\Delta^{d-2}$, we have
\begin{equation}\label{eq:H-positive-definite}
Q_H(\mathbf a,\mathbf a;\mathbf u)
=\frac{\displaystyle\sum_{j=1}^r\ell_j(\mathbf a)^2
\prod_{i=3}^d\ell_j(\mathbf u_i)}
{\displaystyle\sum_{j=1}^r\ell_j(\mathbf v)^d}>0.
\end{equation}
The empty product is interpreted as $1$ when $d=2$.
Each $\ell_j$ is strictly positive on $P\setminus\{\mathbf0\}$, and linear independence ensures that $\ell_j(\mathbf a)\ne0$ for at least one $j$.

\begin{proof}[Proof of Theorem~\ref{thm:lor-main}]
The cases $d=0$, $r=1$, and $d=1$ have already been proved.
Assume $r,d\ge2$ and use the reduced setting above.
In Theorem~\ref{thm:main-criterion}, take
\begin{equation}\label{eq:lorentzian-matrix-data}
C=\Psi(\mathcal C),\qquad
Y=\Psi(\mathcal Y_L),\qquad
B=I_N,\qquad
\mathbf c=J(f_*-H),
\end{equation}
and take $D$ as in~\eqref{eq:cone-diagonal}.
The sets $C$ and $Y$ are closed, $C$ is bounded, and $C\cap Y\ne\varnothing$ because $f_*$ belongs to $\mathcal L_{K,d,\mathbf v}$.
By Lemma~\ref{lem:common-contraction}, the first hypothesis of the matrix criterion holds.

Fix $T>0$, a continuous function $\lambda:[0,T]\to(0,\infty)$, and an arbitrary $f_0\in\mathcal Y_L$.
In polynomial coordinates, the controlled equation is
\[
\frac{d}{dt}f(t)=Lf(t)+\lambda(t)(f(t)-H),
\qquad f(0)=f_0.
\]
All derivatives and integrals of polynomial-valued functions are taken coefficientwise.
The equation preserves $\mathcal A$, because $(Lf)(\mathbf v)=0$ and $(f-H)(\mathbf v)=0$ for $f\in\mathcal A$.
Set $\Lambda(t)=\int_0^t\lambda(s)\,ds$.
By differentiating and using~\eqref{eq:cone-generator}, we have
\[
\frac{d}{dt}\bigl(e^{-\Lambda(t)}T(-t,f(t))\bigr)
=-\lambda(t)e^{-\Lambda(t)}T(-t,H).
\]
Integrating from $0$ to $T$ and applying $e^{\Lambda(T)}T(T,\cdot)$, we obtain
\begin{equation}\label{eq:lorentzian-controlled-solution}
f(T)=e^{\Lambda(T)}T(T,f_0)
-\int_0^T\lambda(s)e^{\Lambda(T)-\Lambda(s)}T(T-s,H)\,ds.
\end{equation}
Replacing $T$ by $t$ in this formula defines a solution on $[0,T]$, as follows by differentiation.
Since every solution satisfies this formula, the solution is unique.

Fix $\mathbf u\in\Delta^{d-2}$.
By Lemma~\ref{lem:lorentzian-intersection}, the mixed form of $e^{\Lambda(T)}T(T,f_0)$ has at most one positive eigenvalue.
By~\eqref{eq:mixed-pullback} and~\eqref{eq:H-positive-definite}, the mixed form of $T(T-s,H)$ is positive definite for every $s\in[0,T]$.
Since $\lambda(s)e^{\Lambda(T)-\Lambda(s)}>0$, the mixed form of the integral in~\eqref{eq:lorentzian-controlled-solution} is positive definite.
Indeed, evaluation on any nonzero $\mathbf a$ is the integral of a continuous strictly positive function on $[0,T]$.
By Proposition~\ref{prop:inertia-facts}, the first form is nonpositive on a subspace of dimension $r-1$.
Subtracting the positive definite form makes the restriction to that subspace negative definite.
By the same proposition, $\theta_2(f(T),\mathbf u)<0$.
This holds for every $\mathbf u\in\Delta^{d-2}$, so~\eqref{eq:lorentzian-interior-test} yields $f(T)\in\operatorname{int}_{\mathcal A}\mathcal Y_L$.

For $\mathbf u(t)=\Psi(f(t))$, the term $Lf(t)$ has coordinate vector $-D\mathbf u(t)$ by~\eqref{eq:diagonal-coordinates}, and
$J(f(t)-H)=\mathbf u(t)+J(f_*-H)$.
Thus the polynomial equation is exactly the matrix equation with the data in~\eqref{eq:lorentzian-matrix-data}.
Since $f_0$ was arbitrary in $\mathcal Y_L$, the second hypothesis holds for every initial point of $Y$.
By Theorem~\ref{thm:main-criterion} and Lemma~\ref{lem:lorentzian-intersection}, the space $\mathcal L_{K,d,\mathbf v}$ in the reduced setting is homeomorphic to $\overline B^N$, where $N=\binom{r+d-1}{d}-1$.
By Corollary~\ref{cor:quotient}, the same conclusion holds in the original space $\mathbb R^n$.
\end{proof}

\subsection{Stable polynomials}
\label{subsec:cone-stable}

We continue to assume that $r\ge2$ and $d\ge2$.
An \emph{entire function} is a function holomorphic on $\mathbb C$, and a \emph{real entire function} is an entire function with real Taylor coefficients at zero.
A nonzero real entire function belongs to the \emph{Laguerre--P\'olya class} if it is a limit of real-rooted real polynomials, uniformly on every compact subset of $\mathbb C$.
For a real entire function $\varphi(z)=\sum_{j\ge0}a_jz^j$, define the constant-coefficient differential operator $\varphi(\partial_z)$ on $\mathbb R[z]$ by
\[
\varphi(\partial_z)p
=
\sum_{j=0}^{\deg p}a_j\frac{d^jp}{dz^j}.
\]
We also set $\varphi(\partial_z)0=0$.

\begin{theorem}[P\'olya--Benz; {\cite[Theorem~1.2]{ABH}}]
\label{thm:polya-benz}
Let $\varphi$ belong to the Laguerre--P\'olya class.
For every real-rooted real polynomial $p$, the polynomial $\varphi(\partial_z)p$ is real-rooted or identically zero.
\end{theorem}

For $a\in\mathbb R$ and nonzero $p\in\mathbb R[z]$, define
\[
e^{-a\partial_z^2}p
=
\sum_{j=0}^{\lfloor\deg p/2\rfloor}
\frac{(-a)^j}{j!}\frac{d^{2j}p}{dz^{2j}}.
\]
The operator sends the zero polynomial to zero.
The expression is a finite sum, and the operator preserves the degree and leading coefficient of every nonzero polynomial.

\begin{corollary}
\label{cor:heat-real-rooted}
Let $p\in\mathbb R[z]$ be nonzero and real-rooted.
For $a\ge0$, the polynomial $e^{-a\partial_z^2}p$ is real-rooted.
If $a>0$, every root of this polynomial is simple.
\end{corollary}

\begin{proof}
For $a=0$, the operator is the identity.
Suppose that $a>0$.
For each positive integer $k$, the polynomial $(1-az^2/k)^k$ is real-rooted, and these polynomials converge to $e^{-az^2}$ uniformly on compact subsets of $\mathbb C$ as $k\to\infty$.
Hence $e^{-az^2}$ belongs to the Laguerre--P\'olya class.
By Theorem~\ref{thm:polya-benz}, the polynomial $e^{-a\partial_z^2}p$ is real-rooted or identically zero.
Since the leading coefficient is unchanged, the polynomial is nonzero.

For $s\ge0$, write $p_s=e^{-s\partial_z^2}p$.
Suppose that $p_a$ has a root $\xi\in\mathbb R$ of multiplicity $k\ge2$.
Write
\[
p_a(\xi+z)=\sum_{j=k}^{\deg p} c_jz^j,
\qquad c_k\ne0.
\]
For $0<h<a$, by the finite exponential identity, we have $p_{a-h}=e^{h\partial_z^2}p_a$.
Since differentiation commutes with translation, we obtain
\[
\begin{aligned}
h^{-k/2}p_{a-h}(\xi+\sqrt h\,w)
&=
\sum_{j=k}^{\deg p} c_jh^{(j-k)/2}
\sum_{\nu=0}^{\lfloor j/2\rfloor}
\frac{j!}{\nu!(j-2\nu)!}w^{j-2\nu}.
\end{aligned}
\]
Consequently, as $h\to0^+$, these polynomials converge in coefficients to
\[
c_k\sum_{\nu=0}^{\lfloor k/2\rfloor}
\frac{k!}{\nu!(k-2\nu)!}w^{k-2\nu}
=
c_ke^{\partial_w^2}w^k.
\]
The polynomial on the left before taking the limit is real-rooted, because $p_{a-h}$ is real-rooted and $\xi$ and $\sqrt h$ are real.
The limit is nonzero.
By Theorem~\ref{thm:real-rooted-limits}, the polynomial $e^{\partial_w^2}w^k$ is real-rooted.

This polynomial is monic, with coefficient zero at $w^{k-1}$ and coefficient $k(k-1)>0$ at $w^{k-2}$.
If $\xi_1,\ldots,\xi_k$ were all real roots, counted with multiplicity, by comparison of coefficients, we would have
\[
\sum_{j=1}^k\xi_j=0,
\qquad
\sum_{1\le i<j\le k}\xi_i\xi_j=k(k-1)>0.
\]
However,
\[
\sum_{1\le i<j\le k}\xi_i\xi_j
=
\frac12\left(\left(\sum_{j=1}^k\xi_j\right)^2-
\sum_{j=1}^k\xi_j^2\right)
=
-\frac12\sum_{j=1}^k\xi_j^2\le0,
\]
a contradiction.
Therefore every root of $p_a$ is simple.
\end{proof}

Define
\begin{equation}\label{eq:hyperbolic-structural-set}
\mathcal Y_{\mathrm{hyp}}
=
\{f\in\mathcal A:
 f\text{ is hyperbolic with respect to }\mathbf v\}.
\end{equation}

\begin{lemma}
\label{lem:stable-interior}
The set $\mathcal Y_{\mathrm{hyp}}$ is closed in $\mathcal A$, and
\[
\mathcal C\cap\mathcal Y_{\mathrm{hyp}}
=
\mathcal S_{K,d,\mathbf v}.
\]
Every polynomial in $\mathcal A$ that is strictly hyperbolic with respect to $\mathbf v$ belongs to $\operatorname{int}_{\mathcal A}\mathcal Y_{\mathrm{hyp}}$.
For every $t\in\mathbb R$, the map $f\mapsto T(t,f)$ preserves both hyperbolicity and strict hyperbolicity with respect to $\mathbf v$.
\end{lemma}

\begin{proof}
Let $(f_j)_{j\ge1}$ be a sequence in $\mathcal Y_{\mathrm{hyp}}$ converging to $f\in\mathcal A$.
For every $\mathbf x\in\mathbb R^r$, the polynomials $f_j(\mathbf x+z\mathbf v)$ converge in coefficients to $f(\mathbf x+z\mathbf v)$.
All these polynomials are monic of degree $d$, since $f_j(\mathbf v)=f(\mathbf v)=1$.
By Theorem~\ref{thm:real-rooted-limits}, the polynomial $f(\mathbf x+z\mathbf v)$ is real-rooted.
Hence $f\in\mathcal Y_{\mathrm{hyp}}$, so $\mathcal Y_{\mathrm{hyp}}$ is closed.

By Lemma~\ref{lem:common-compact} and Corollary~\ref{cor:stable-hyperbolic}, we have
$\mathcal S_{K,d,\mathbf v}\subseteq\mathcal C\cap\mathcal Y_{\mathrm{hyp}}$.
Conversely, suppose that $f\in\mathcal C\cap\mathcal Y_{\mathrm{hyp}}$.
By Lemma~\ref{lem:common-compact}, we have $f(\mathbf x)>0$ for every $\mathbf x\in K$.
The set $K$ is connected and contains $\mathbf v$.
Hence $K$ is contained in the connected component of $\{\mathbf x:f(\mathbf x)\ne0\}$ containing $\mathbf v$, namely $\Gamma_{\mathbf v}(f)$.
By Corollary~\ref{cor:stable-hyperbolic}, we have $f\in\mathcal S_{K,d,\mathbf v}$.
This proves the asserted equality.

Let $f\in\mathcal A$ be strictly hyperbolic with respect to $\mathbf v$, and let
\[
\Sigma=\{\mathbf y\in W:\|\mathbf y\|=1\}.
\]
For every $\mathbf y\in\Sigma$, the monic polynomial $f(\mathbf y+z\mathbf v)$ has $d$ distinct real roots.
The set of monic real polynomials of degree $d$ with distinct real roots is open in the coefficient space.
Indeed, for each such polynomial, choose disjoint intervals about the $d$ roots such that the polynomial has nonzero values of opposite signs at the two endpoints of each interval.
These signs remain unchanged under a sufficiently small coefficient perturbation.
By the intermediate value theorem, the perturbed polynomial has a root in each of the $d$ intervals.
Since the degree is $d$, all roots are real and simple.

The map that sends $(g,\mathbf y)\in\mathcal A\times\Sigma$ to the coefficient vector of $g(\mathbf y+z\mathbf v)$ is continuous.
For each $\mathbf y\in\Sigma$, there are neighborhoods $V_{\mathbf y}$ of $f$ in $\mathcal A$ and $O_{\mathbf y}$ of $\mathbf y$ in $\Sigma$ such that $g(\mathbf w+z\mathbf v)$ has $d$ distinct real roots whenever $g\in V_{\mathbf y}$ and $\mathbf w\in O_{\mathbf y}$.
Since $\Sigma$ is compact, finitely many of the sets $O_{\mathbf y}$ cover $\Sigma$.
The intersection of the corresponding sets $V_{\mathbf y}$ is a neighborhood $V$ of $f$ in $\mathcal A$ on which the restrictions are real-rooted for all $\mathbf y\in\Sigma$.
For $g\in V$ and $\mathbf y\in W\setminus\{\mathbf0\}$, by homogeneity, we have
\[
g(\mathbf y+z\mathbf v)
=
\|\mathbf y\|^d
 g\left(\frac{\mathbf y}{\|\mathbf y\|}
 +\frac{z}{\|\mathbf y\|}\mathbf v\right).
\]
Hence $g(\mathbf y+z\mathbf v)$ is real-rooted for all nonzero $\mathbf y\in W$.
For $\mathbf y=\mathbf0$, this restriction is $z^d$.
Every vector $\mathbf x\in\mathbb R^r$ has the form $\mathbf y+s\mathbf v$ with $\mathbf y\in W$ and $s\in\mathbb R$.
Translation of the variable $z$ therefore proves that $g$ is hyperbolic with respect to $\mathbf v$.
Thus $V\subseteq\mathcal Y_{\mathrm{hyp}}$ and $f\in\operatorname{int}_{\mathcal A}\mathcal Y_{\mathrm{hyp}}$.

For every $t\in\mathbb R$, the map $G(t,\cdot)$ is real linear, invertible, and fixes $\mathbf v$.
Consequently,
\[
T(t,f)(\mathbf x+z\mathbf v)
=
f(G(t,\mathbf x)+z\mathbf v).
\]
Hyperbolicity with respect to $\mathbf v$ is therefore preserved by $T(t,\cdot)$.
Moreover, $G(t,\mathbf x)\in\mathbb R\mathbf v$ if and only if $\mathbf x\in\mathbb R\mathbf v$.
The same identity proves preservation of strict hyperbolicity.
\end{proof}

Define $q\in\mathrm H_2(\mathbb R^r)$ and a linear operator $\mathcal H:\mathrm H_d(\mathbb R^r)\to\mathrm H_d(\mathbb R^r)$ by
\begin{equation}\label{eq:heat-operator}
q(\mathbf x)=\|\mathbf x-\ell(\mathbf x)\mathbf v\|^2,
\qquad
\mathcal Hf=qD_{\mathbf v}^{\,2}f.
\end{equation}
The vector $\mathbf x-\ell(\mathbf x)\mathbf v$ belongs to $W$.
Since $\ell(\mathbf v)=1$, we have $q(\mathbf x+s\mathbf v)=q(\mathbf x)$ for every $s\in\mathbb R$, so $D_{\mathbf v}q=0$.
Also, $q(\mathbf x)=0$ if and only if $\mathbf x\in\mathbb R\mathbf v$.
For every integer $j\ge0$, we have
\[
\mathcal H^jf=q^jD_{\mathbf v}^{\,2j}f.
\]
In particular, $\mathcal H^{\lfloor d/2\rfloor+1}=0$.
An operator for which a positive integer power is zero is called \emph{nilpotent}.
For $a\in\mathbb R$, define
\begin{equation}\label{eq:heat-finite}
e^{-a\mathcal H}f
=
\sum_{j=0}^{\lfloor d/2\rfloor}
\frac{(-a)^j}{j!}q^jD_{\mathbf v}^{\,2j}f.
\end{equation}
Since $q(\mathbf v)=0$, we have $(\mathcal Hf)(\mathbf v)=0$ and $(e^{-a\mathcal H}f)(\mathbf v)=f(\mathbf v)$.
Hence $\mathcal H$ maps $\mathrm H_d(\mathbb R^r)$ into $E$, and $e^{-a\mathcal H}$ preserves $\mathcal A$.

\begin{corollary}
\label{cor:heat-hyperbolic}
For $f\in\mathcal Y_{\mathrm{hyp}}$ and $a\ge0$, we have $e^{-a\mathcal H}f\in\mathcal Y_{\mathrm{hyp}}$.
If $a>0$, the polynomial $e^{-a\mathcal H}f$ is strictly hyperbolic with respect to $\mathbf v$.
\end{corollary}

\begin{proof}
Fix $\mathbf y\in W$ and write $p(z)=f(\mathbf y+z\mathbf v)$.
The polynomial $p$ is monic of degree $d$ and real-rooted.
Since $q(\mathbf y+z\mathbf v)=\|\mathbf y\|^2$, by~\eqref{eq:heat-finite}, we have
\[
(e^{-a\mathcal H}f)(\mathbf y+z\mathbf v)
=
e^{-a\|\mathbf y\|^2\partial_z^2}p(z).
\]
By Corollary~\ref{cor:heat-real-rooted}, this polynomial is real-rooted for $a\ge0$ and has $d$ distinct real roots when $a>0$ and $\mathbf y\ne\mathbf0$.
Every vector in $\mathbb R^r$ has the form $\mathbf y+s\mathbf v$ with $\mathbf y\in W$.
Translation of the univariate argument proves both assertions.
\end{proof}

\begin{proof}[Proof of Theorem~\ref{thm:stable-main}]
The cases $d=0$, $r=1$, and $d=1$ have already been proved.
By Corollary~\ref{cor:quotient}, it remains to consider a pointed cone in $\mathbb R^r$, with $r\ge2$ and $d\ge2$.
In the affine coordinates $\Psi$, set
\begin{equation}\label{eq:stable-matrix-data}
C=\Psi(\mathcal C),\qquad
Y=\Psi(\mathcal Y_{\mathrm{hyp}}),\qquad
B=-J\circ(\mathcal H|_E)\circ J^{-1},\qquad
\mathbf c=-J(\mathcal Hf_*).
\end{equation}
We identify the linear map $B:\mathbb R^N\to\mathbb R^N$ with the matrix in the standard basis.
The diagonal matrix $D$ is the matrix from Lemma~\ref{lem:common-contraction}.
By Lemmas~\ref{lem:common-compact} and~\ref{lem:stable-interior}, the sets $C$ and $Y$ are closed, and $C$ is bounded.
The polynomial $f_*=\ell^d$ belongs to $\mathcal S_{K,d,\mathbf v}$, so $C\cap Y$ is nonempty.
By Lemma~\ref{lem:common-contraction}, the flow with $\lambda\equiv0$ maps $C$ into $\operatorname{int}C$ at every positive time.

Fix $T>0$, a continuous function $\lambda:[0,T]\to(0,\infty)$, and an arbitrary $f_0\in\mathcal Y_{\mathrm{hyp}}$.
Under the coordinate map $\Psi$, the matrix equation in Theorem~\ref{thm:main-criterion} becomes
\begin{equation}\label{eq:stable-polynomial-equation}
\frac{d}{dt}f(t)
=
Lf(t)-\lambda(t)\mathcal Hf(t),
\qquad f(0)=f_0.
\end{equation}
Since $Lf\in E$ and $\mathcal Hf\in E$, this equation preserves $\mathcal A$.
Indeed, $J(Lf)=-D\Psi(f)$ and
\[
J(-\mathcal Hf)
=
-J\bigl(\mathcal H(f-f_*)\bigr)-J(\mathcal Hf_*)
=
B\Psi(f)+\mathbf c.
\]
By~\cite[Theorem~2.2 and Corollary~2.6]{Teschl}, the equation has a unique solution on $[0,T]$, since it is a linear differential equation on a finite-dimensional coefficient space with continuous coefficients.

For every $f\in\mathrm H_d(\mathbb R^r)$, by direct differentiation, we have
\[
\begin{aligned}
L\mathcal Hf
&=\ell qD_{\mathbf v}^{\,3}f-dqD_{\mathbf v}^{\,2}f,\\
\mathcal HLf
&=q\left(2D_{\mathbf v}^{\,2}f
 +\ell D_{\mathbf v}^{\,3}f-dD_{\mathbf v}^{\,2}f\right).
\end{aligned}
\]
Thus $L\mathcal H-\mathcal HL=-2\mathcal H$.
Differentiating $T(t,\mathcal H T(-t,f))$ with respect to $t$, and using the value at $t=0$, we obtain
\begin{equation}\label{eq:heat-conjugation}
T(t,\mathcal H T(-t,f))=e^{-2t}\mathcal Hf
\qquad(t\in\mathbb R).
\end{equation}
By~\eqref{eq:stable-polynomial-equation} and~\eqref{eq:heat-conjugation}, we have
\[
\frac{d}{dt}T(-t,f(t))
=
-\lambda(t)e^{2t}\mathcal H T(-t,f(t)).
\]
Define
\[
a(t)=\int_0^t\lambda(s)e^{2s}\,ds
\qquad(0\le t\le T).
\]
The function $e^{-a(t)\mathcal H}f_0$ solves the preceding differential equation with initial value $f_0$.
By uniqueness, we have
\begin{equation}\label{eq:stable-controlled-solution}
f(t)=T(t,e^{-a(t)\mathcal H}f_0)
\qquad(0\le t\le T).
\end{equation}
Since $\lambda$ is positive, we have $a(T)>0$.
By Corollary~\ref{cor:heat-hyperbolic}, the polynomial $e^{-a(T)\mathcal H}f_0$ is strictly hyperbolic with respect to $\mathbf v$.
By Lemma~\ref{lem:stable-interior}, the polynomial $f(T)$ is also strictly hyperbolic and belongs to $\operatorname{int}_{\mathcal A}\mathcal Y_{\mathrm{hyp}}$.
Since $\Psi$ is a homeomorphism of the ambient affine spaces, we have $\Psi(f(T))\in\operatorname{int}Y$.
This verifies the second hypothesis of Theorem~\ref{thm:main-criterion} for every initial point in $Y$.

By Theorem~\ref{thm:main-criterion}, the set $C\cap Y$ is homeomorphic to $\overline B^N$.
By Lemma~\ref{lem:stable-interior}, we have
\[
C\cap Y=\Psi(\mathcal S_{K,d,\mathbf v}).
\]
By the dimension formula for $\mathcal A$ and the homeomorphism in Corollary~\ref{cor:quotient}, the original normalized stable space is homeomorphic to a closed Euclidean ball of dimension $\binom{r+d-1}{d}-1$.
When this dimension is positive, the resulting homeomorphism maps the boundary relative to $\mathcal A_{K,d,\mathbf v}$ onto the Euclidean unit sphere.
\end{proof}

\section{Multiaffine real stable polynomials}\label{sec:ma}

\subsection{The coefficient cone and preserving operators}

Write $[n]=\{1,\ldots,n\}$, and let $\binom{[n]}d$ denote the set of
$d$-element subsets of $[n]$. For $I\subseteq[n]$, put
$x^I=\prod_{i\in I}x_i$, with an empty product equal to $1$.
Define
\begin{align*}
V_{n,d}^{\ma}=\left\{\sum_{I\in\binom{[n]}d}a_Ix^I:a_I\in\mathbb R\right\}.
\end{align*}
This is the vector space of homogeneous multiaffine polynomials of
degree $d$, of dimension $m=\binom nd$.
For $f=\sum_I a_Ix^I$ and $g=\sum_I b_Ix^I$ in
$V_{n,d}^{\mathrm{ma}}$, define the coefficient inner product by
\begin{align*}
\langle f,g\rangle
:=\sum_{I\in\binom{[n]}d}a_Ib_I.
\end{align*}
 In particular, for each $I\in\binom{[n]}d$, let $e_I=x^I$.
 These monomials form an orthonormal basis of
 $V_{n,d}^{\mathrm{ma}}$, since
 \begin{align*}
 \langle e_I,e_J\rangle=
 \begin{cases}
 	1,& I=J,\\
 	0,& I\ne J,
 \end{cases}
 \qquad I,J\in\binom{[n]}d.
 \end{align*}

All operators below are linear maps on $V_{n,d}^{\mathrm{ma}}$,
represented in the monomial basis
$(e_I)_{I\in\binom{[n]}d}$.
For such an operator $T$, write
\begin{align*}
T(e_J)=\sum_{I\in\binom{[n]}d}T_{IJ}e_I,
\qquad J\in\binom{[n]}d.
\end{align*}
Then $(T_{IJ})_{I,J}$ is the matrix of $T$ in this basis.

We first prove Theorem~\ref{thm:ma-main} in the cases $d=0$ and $d=n$.

Suppose first that $d=0$. Every polynomial
$f\in\mathcal M_{n,0}$ is constant, so $f=c$ for some
$c\in\mathbb R$. The normalization $f(\mathbf1)=1$ gives
$c=1$. Conversely, the constant polynomial $1$ is real stable
and satisfies the normalization. Therefore,
$\mathcal M_{n,0}=\{1\}.$

Next suppose that $d=n$. Since $f$ is homogeneous of degree $n$
and multiaffine, each monomial has total degree $n$ and contains
each variable at most once. Hence each monomial must contain
all $n$ variables exactly once, so $f=c\,x_1\cdots x_n$.
The normalization $f(\mathbf1)=1$ then gives $c=1$.
Conversely, $x_1\cdots x_n$ is homogeneous of degree $n$,
multiaffine, and equal to $1$ at $\mathbf1$. It is also real
stable: if $\operatorname{Im}z_i>0$ for every $i$, then every
$z_i$ is nonzero, and hence $z_1\cdots z_n\ne0$.
Consequently, $\mathcal M_{n,n}=\{x_1\cdots x_n\}.$

Both spaces consist of a single point and are therefore
homeomorphic to a zero-dimensional closed ball.
Since $\binom n0-1=\binom nn-1=0$, this proves
Theorem~\ref{thm:ma-main} for $d=0$ and $d=n$.
We henceforth assume $1\le d\le n-1$.

We first justify the coefficient condition used in this section.
\begin{theorem}[{\cite[Theorem~6.1]{COSW}}]
\label{thm:same-phase}
If a homogeneous complex polynomial is nonvanishing whenever every
variable has positive real part, all its nonzero coefficients have the
same phase: the ratio of any two such coefficients is a positive real
number.
\end{theorem}

\begin{corollary}\label{cor:normalized-coefficients}
	Let $f$ be a homogeneous real stable polynomial.
	Then all its nonzero coefficients have the same sign.
	In addition, if $f(\mathbf1)=1$, then all its coefficients
	are nonnegative.
\end{corollary}

\begin{proof}
	Let $d$ be the degree of $f$, and take
	$z=(z_1,\ldots,z_n)\in\mathbb C^n$ with
	$\operatorname{Re}z_j>0$ for every $j$.
	Since
	\begin{align*}
		\operatorname{Im}(\mathrm{i}z_j)
		=\operatorname{Re}z_j>0,
	\end{align*}
	every coordinate of $\mathrm{i}z$ lies in the upper half-plane.
	Real stability therefore gives $f(\mathrm{i}z)\ne0$.
	On the other hand, homogeneity gives
	\begin{align*}
		f(\mathrm{i}z)=\mathrm{i}^d f(z).
	\end{align*}
	Because $\mathrm{i}^d\ne0$, we obtain $f(z)\ne0$.
	Thus $f$ is nonvanishing whenever all variables have positive
	real part. By Theorem~\ref{thm:same-phase}, the ratio of any
	two nonzero coefficients of $f$ is a positive real number.
	Since these coefficients are real, they must have the same sign.
	
	Now assume $f(\mathbf1)=1$. Each monomial equals $1$ at
	$\mathbf1$, so the sum of all coefficients is $1$.
	The nonzero coefficients therefore cannot all be negative.
	Since they have the same sign, they must all be positive.
	The remaining coefficients are zero, so all coefficients
	of $f$ are nonnegative.
\end{proof}

To study $\mathcal M_{n,d}$ using linear operators, we first
consider the following larger set, without imposing the
normalization $f(\mathbf1)=1$:
\begin{align}\label{eq:ma-cone}
	\mathscr C_{n,d}
	:=\{0\}\cup
	\{f\in V_{n,d}^{\ma}:
	f\text{ is real stable and has nonnegative coefficients}\}.
\end{align}

By Corollary~\ref{cor:normalized-coefficients}, every
$f\in\mathcal M_{n,d}$ has nonnegative coefficients.
Hence $f\in\mathscr C_{n,d}$ and $f(\mathbf1)=1$.
Conversely, if $f\in\mathscr C_{n,d}$ and $f(\mathbf1)=1$,
then $f\ne0$ and satisfies all the defining conditions of
$\mathcal M_{n,d}$. Therefore,
\begin{align*}
	\mathcal M_{n,d}
	=\{f\in\mathscr C_{n,d}:f(\mathbf1)=1\}.
\end{align*}

We next verify that $\mathscr C_{n,d}$ is closed.
By~\cite[Proposition~3.1(3)]{BBL}, a coefficient limit
of real stable polynomials is either real stable or zero.
Since nonnegative coefficients remain nonnegative under
limits, every limit in $V_{n,d}^{\ma}$ of polynomials in
$\mathscr C_{n,d}$ still belongs to $\mathscr C_{n,d}$.

Finally, let $f\in\mathscr C_{n,d}$ and $t\ge0$.
If $f=0$ or $t=0$, then $tf=0\in\mathscr C_{n,d}$.
Otherwise, $f$ is real stable with nonnegative coefficients
and $t>0$. Multiplication by $t$ preserves both properties,
so $tf\in\mathscr C_{n,d}$.
Thus $\mathscr C_{n,d}$ is closed under multiplication
by nonnegative real numbers and is therefore a cone.
We do not assume that it is convex.

By~\cite[Proposition~3.1(2)]{BBL}, every partial derivative
of a real stable polynomial is either real stable or zero.
We use this fact below to prove that the following operators
preserve $\mathscr C_{n,d}$.

Let $\Id$ denote the identity operator on $V_{n,d}^{\ma}$.
For $i\in[n]$, define $$N_i=x_i\partial_i.$$
For distinct $i,j\in[n]$, define
$$E_{ij}=x_i\partial_j(\Id-N_i).$$
For $I\in\binom{[n]}d$, we have
\begin{align*}
	N_ix^I=
	\begin{cases}
		x^I,i\in I,\\
		0,i\notin I,
	\end{cases}\qquad
	E_{ij}x^I=
	\begin{cases}
		x^{(I\setminus\{j\})\cup\{i\}},j\in I,\ i\notin I,\\
		0,\text{otherwise}.
	\end{cases}
\end{align*}
Consequently, $N_i^2=N_i$, $E_{ij}^2=0$, and
$\sum_{i=1}^nN_i=d\,\Id$.

For a linear operator $A$ on $V_{n,d}^{\ma}$ and $t\in\mathbb R$,
define its exponential by
\begin{align*}
	e^{tA}:=\sum_{k=0}^{\infty}\frac{t^k}{k!}A^k,
	\qquad A^0=\Id.
\end{align*}

\begin{lemma}\label{lem:elementary-preservers}
	Let $f\in\mathscr C_{n,d}$.
	For every $i\in[n]$ and $s\in\mathbb R$, we have
	\begin{align*}
		e^{sN_i}f\in\mathscr C_{n,d}.
	\end{align*}
	For distinct $i,j\in[n]$ and every $t\ge0$, we have
	\begin{align*}
		e^{tE_{ij}}f=f+tE_{ij}f\in\mathscr C_{n,d}.
	\end{align*}

\end{lemma}

\begin{proof}
	Let $f\in\mathscr C_{n,d}$. If $f=0$, both operators send
	$f$ to zero, so we may assume $f\ne0$.
	
	First, by $N_i^2=N_i$, we have
	\begin{align*}
		e^{sN_i}f
		=f+(e^s-1)N_if
		=f(x_1,\ldots,e^sx_i,\ldots,x_n).
	\end{align*}
	Since $e^s>0$, this substitution preserves the upper half-plane
	and multiplies coefficients by positive numbers.
	Hence $e^{sN_i}f$ is real stable and has nonnegative coefficients,
	so $e^{sN_i}f\in\mathscr C_{n,d}$.
	
	Next, let $t\ge0$ and $i\ne j$, and set
	$
		p(x)=f(x_1,\ldots,x_j+tx_i,\ldots,x_n).
	$
	If all variables lie in $\mathbb H$, then
	$\operatorname{Im}(x_j+tx_i)>0$. Thus $p$ is real stable.
	Since $f$ is multiaffine, its degree in $x_j$ is at most one.
	By the definitions of $N_i$ and $E_{ij}$, we therefore have
	\begin{align*}
		p=f+tx_i\partial_jf
		=f+tE_{ij}f+tx_i^2\partial_i\partial_jf.
	\end{align*}
	The polynomial $\partial_i\partial_jf$ does not depend on $x_i$,
	whereas $f+tE_{ij}f$ has degree at most one in $x_i$.
	Consequently, by deleting the quadratic terms in $x_i$ from $p$,
	we obtain
	$
		r=f+tE_{ij}f.
	$
	
	We now verify that $r$ is stable or zero.
	Write $z=x_i$ and
	$
		p(z)=a+bz+cz^2,
	$
	where $a,b,c$ are polynomials in the other variables.
	Since $z\mapsto-1/z$ maps $\mathbb H$ into itself and
	$z\ne0$ on $\mathbb H$, the polynomial
	$
		q(z)=z^2p(-1/z)=az^2-bz+c
	$
	is stable. By~\cite[Proposition~3.1(2)]{BBL},
	$q'$ is stable or zero.
	Since $z\mapsto-1/(2z)$ also maps $\mathbb H$ into itself,
	we obtain
	\begin{align*}
		r(z)=a+bz=(-z)q'\left(-\frac1{2z}\right),
	\end{align*}
	which is stable or zero.
	
	By $E_{ij}^2=0$, we have
	$
		(\Id-tE_{ij})(\Id+tE_{ij})=\Id.
	$
	Since $f\ne0$, this identity implies
	$r=(\Id+tE_{ij})f\ne0$. Hence $r$ is real stable.
	Moreover, $p$ has nonnegative coefficients because $f$ does
	and $t\ge0$. Deleting its quadratic terms preserves coefficient
	nonnegativity. Thus $r\in\mathscr C_{n,d}$.
	Finally, $E_{ij}^2=0$ implies
	$e^{tE_{ij}}f=(\Id+tE_{ij})f=r$, proving the second assertion.
\end{proof}

For real numbers $h_i$ and positive numbers $w_{ij}=w_{ji}$, define
\begin{equation}\label{eq:exchange-matrix}
 \mathsf M=\sum_{i=1}^nh_iN_i+\sum_{1\le i<j\le n}w_{ij}(E_{ij}+E_{ji}).
\end{equation}
The matrix $\mathsf M$ is real symmetric in the monomial basis,
and its off-diagonal entries are nonnegative.
Consider the graph whose vertices are the degree-$d$ monomials,
with two vertices joined when one monomial is obtained from the
other by a single exchange.
Since every $w_{ij}>0$, each edge corresponds to a positive
off-diagonal entry of $\mathsf M$.
This graph is connected: for any $I,J\in\binom{[n]}d$,
we can transform $I$ into $J$ by successively replacing elements
of $I\setminus J$ with elements of $J\setminus I$.

\begin{theorem}[{\cite[Lie--Trotter product formula]{Trotter}}]\label{thm:trotter}
For finitely many linear maps $A_1,\ldots,A_b$ on a finite-dimensional
space and $t\ge0$,
\begin{equation}\label{eq:product-formula}
 e^{t\sum_jA_j}
   =\lim_{k\to\infty}(e^{tA_1/k}\cdots e^{tA_b/k})^k
\end{equation}
in operator norm.
\end{theorem}

\begin{lemma}\label{lem:controlled-preservers}
	Let $T>0$ and let $u_1,\ldots,u_n:[0,T]\to\R$ be continuous.
	For every $y_0\in\mathscr C_{n,d}$, the solution of
	\begin{align}\label{eq:controlled-ode}
		\dot y(s)
		&=\left(\mathsf M+\sum_{i=1}^n u_i(s)N_i\right)y(s),
		&
		y(0)&=y_0
	\end{align}
	satisfies $y(s)\in\mathscr C_{n,d}$ for every $s\in[0,T]$.
	In particular, taking $u_i=0$ for all $i$, we obtain
	$e^{t\mathsf M}f\in\mathscr C_{n,d}$ for every
	$f\in\mathscr C_{n,d}$ and $t\ge0$.
\end{lemma}

\begin{proof}
	First suppose that each control is constant, say $u_i(s)=v_i$.
	Then the coefficient matrix is
	\begin{align*}
		A_v
		=\mathsf M+\sum_{i=1}^n v_iN_i
		=\sum_{i=1}^n(h_i+v_i)N_i
		+\sum_{1\le i<j\le n}w_{ij}(E_{ij}+E_{ji}).
	\end{align*}
	By Lemma~\ref{lem:elementary-preservers},
	$e^{\tau(h_i+v_i)N_i}$ preserves $\mathscr C_{n,d}$
	for every $\tau\ge0$, since $h_i+v_i$ may be any real number.
	The same lemma applies to $e^{\tau w_{ij}E_{ij}}$ and
	$e^{\tau w_{ij}E_{ji}}$, since $w_{ij}>0$.
	Thus every product in \eqref{eq:product-formula} preserves
	$\mathscr C_{n,d}$.
	Since $\mathscr C_{n,d}$ is closed, its limit
	$e^{\tau A_v}$ also preserves $\mathscr C_{n,d}$.
	
Now suppose that $u_1,\ldots,u_n$ are continuous on $[0,T]$.
For each $i$, choose piecewise constant functions $u_i^{(k)}$
converging uniformly to $u_i$ on $[0,T]$, and set
\begin{align*}
	A(s)=\mathsf M+\sum_{i=1}^n u_i(s)N_i,
	\qquad
	A_k(s)=\mathsf M+\sum_{i=1}^n u_i^{(k)}(s)N_i.
\end{align*}
	Let $y_k$ solve $\dot y_k=A_k(s)y_k$ with $y_k(0)=y_0$.
	On each interval where the controls are constant,
	the solution is obtained by applying an exponential
	of the form already considered.
	Starting from $y_0\in\mathscr C_{n,d}$ and proceeding
	through these intervals, we obtain
	$y_k(s)\in\mathscr C_{n,d}$ for every $s\in[0,T]$.
	
	We next verify that $y_k$ converges uniformly to $y$.
	Since $u_i^{(k)}\to u_i$ uniformly, we have
	\begin{align*}
		\delta_k:=\sup_{s\in[0,T]}\|A_k(s)-A(s)\|\longrightarrow0.
	\end{align*}
	Since $A_k$ converges uniformly to $A$ and $y$ is continuous,
	there are constants $L,R$ such that
	$\|A_k(s)\|\le L$ and $\|y(s)\|\le R$ for all $k$ and $s$.
	By the integral equations for $y_k$ and $y$, we have
	\begin{align*}
		y_k(s)-y(s)
		=\int_0^s A_k(r)\bigl(y_k(r)-y(r)\bigr)\,dr
		+\int_0^s\bigl(A_k(r)-A(r)\bigr)y(r)\,dr.
	\end{align*}
	Consequently,
	\begin{align*}
		\|y_k(s)-y(s)\|
		\le L\int_0^s\|y_k(r)-y(r)\|\,dr+TR\delta_k.
	\end{align*}
	By Gr\"onwall's inequality~\cite[Lemma~2.7]{Teschl}, we have 
	\begin{align*}
		\sup_{s\in[0,T]}\|y_k(s)-y(s)\|
		\le TR\delta_k e^{LT}\longrightarrow0.
	\end{align*}
	Since each $y_k(s)$ belongs to the closed set
	$\mathscr C_{n,d}$, we conclude that
	$y(s)\in\mathscr C_{n,d}$ for every $s\in[0,T]$.
	
	Finally, when all controls are zero, the solution is
	$y(s)=e^{s\mathsf M}y_0$.
	Since $T>0$ is arbitrary, this proves the last assertion.
\end{proof}
\subsection{Eigenvalues and eigenvectors of $\mathsf M$}
The \emph{spectral radius} of a matrix is the largest absolute value
of its complex eigenvalues.  A nonnegative matrix is \emph{irreducible}
if its directed graph of positive entries is strongly connected, meaning
that a directed path joins every ordered pair of vertices.

\begin{theorem}[{\cite[Chapter 2]{BermanPlemmons}}]
\label{thm:perron}
An irreducible nonnegative real matrix of order at least two has a positive spectral-radius
eigenvalue of algebraic multiplicity one and strictly positive left and
right eigenvectors for this eigenvalue.
\end{theorem}
An eigenvalue is called \emph{simple} if it has algebraic
multiplicity one, that is, if it is a root of multiplicity one
of the characteristic polynomial.

Since $\mathsf M$ has nonnegative off-diagonal entries,
we can choose $\alpha>0$ sufficiently large so that
$\mathsf Q:=\mathsf M+\alpha\Id$ has nonnegative entries.
Adding $\alpha\Id$ changes only the diagonal entries.
Thus, by symmetry and the connectedness of the exchange graph
described above, $\mathsf Q$ is irreducible.
By Theorem~\ref{thm:perron}, its spectral radius $\rho$
is a simple eigenvalue with a strictly positive eigenvector $p$.
We normalize $p$ so that $\|p\|=1$.

Since $\mathsf M$ is symmetric, all its eigenvalues are real.
Moreover, adding $\alpha\Id$ shifts each eigenvalue by $\alpha$
without changing its multiplicity or eigenvectors.
Therefore, the largest eigenvalue of $\mathsf M$ is
$\lambda_0=\rho-\alpha$, which is simple and satisfies
$\mathsf Mp=\lambda_0p$.
We call $p$ the \emph{Perron vector} of $\mathsf M$.
The following lemma strengthens this conclusion by a suitable
choice of the parameters in \eqref{eq:exchange-matrix}.

\begin{lemma}\label{lem:matrix}
For $1\le d\le n-1$, there are parameters $h_i\in\R$ and $w_{ij}>0$ in
\eqref{eq:exchange-matrix} such that $\mathsf M$ has eigenvalues and
orthonormal eigenvectors
\[
 \lambda_0>\lambda_1>\cdots>\lambda_{m-1},\qquad
 v_0=p,v_1,\ldots,v_{m-1},
\]
with the following properties:
\begin{enumerate}
\item every coordinate of $p$ is strictly positive;
\item if $a\ne b$, $c\ne e$, and $(a,b)\ne(c,e)$, then
$\lambda_a-\lambda_b\ne\lambda_c-\lambda_e$;
\item for each $k\ge1$ there is $i\in[n]$ with
$\ip{v_k}{N_ip}\ne0$.
\end{enumerate}
\end{lemma}
Property (ii) states that all nonzero ordered eigenvalue
differences are distinct.
Lemma~\ref{lem:matrix} is proved in Appendix~\ref{app:matrix}.

\subsection{Endpoint rank and strict regularization}

We state the control-theoretic input before applying it.
A \emph{control} on $[0,t]$ is a continuous function from that interval
to a specified open neighborhood of $0$ in $\R^s$.
The \emph{endpoint map} sends a control to the value at time $t$ of
the solution with a fixed initial value.  Its differential is
\emph{surjective}, or of \emph{full rank}, if its image is the entire
finite-dimensional state space.

\begin{theorem}[{\cite[Section 3.7,
Theorem 7]{Sontag}}]\label{thm:endpoint-input}
For a continuously differentiable control system, if the differential
of the endpoint map with respect to controls along a trajectory is
surjective, then endpoints obtained by controls near the given control
contain a neighborhood of the given endpoint.
\end{theorem}

Let $V$ be an $m$-dimensional real vector space, and let
$\operatorname{End}(V)$ denote the space of all linear maps
from $V$ to itself.
For $M,B\in\operatorname{End}(V)$, define their commutator by
$[M,B]=MB-BM$.
For fixed $M$, define the linear map
$\operatorname{ad}_M:\operatorname{End}(V)\to\operatorname{End}(V)$
by $\operatorname{ad}_M(B)=[M,B]$.
Its iterates satisfy
\begin{align*}
	\operatorname{ad}_M^0(B)=B,
	\qquad
	\operatorname{ad}_M^{k+1}(B)
	=[M,\operatorname{ad}_M^k(B)]
	\quad(k\ge0).
\end{align*}
For $B_1,\ldots,B_s\in\operatorname{End}(V)$, define
\begin{align}\label{eq:commutator-span}
	\mathscr K(M;B_1,\ldots,B_s)
	=\spanop\{\operatorname{ad}_M^k(B_i):
	1\le i\le s,\ 0\le k\le m^2-1\}.
\end{align}

Since $\operatorname{End}(V)$ has dimension $m^2$,
the Cayley--Hamilton theorem
\cite{HornJohnson} implies that
$\operatorname{ad}_M^{m^2}$ is a linear combination of
$\operatorname{ad}_M^0,\ldots,\operatorname{ad}_M^{m^2-1}$.
By induction, the same holds for every higher power.
Thus allowing all $k\ge0$ in \eqref{eq:commutator-span}
does not change the span.

For a linear subspace $\mathscr K\subseteq\operatorname{End}(V)$
and $f\in V$, we write
$
	\mathscr Kf:=\{Bf:B\in\mathscr K\}.
$

\begin{proposition}\label{prop:rank-regularization}
	Let $Z\subseteq V$ be closed, and let $M,B_1,\ldots,B_s$
	be linear maps on $V$.
	Assume that for every $T>0$ there exists $\varepsilon_T>0$
	such that the following holds:
	if $u_1,\ldots,u_s:[0,T]\to\R$ are continuous and
	$|u_i(\tau)|<\varepsilon_T$ for every $i$ and $\tau\in[0,T]$,
	then the solution of
	\begin{align*}
		\dot y(\tau)
		=\left(M+\sum_{i=1}^s u_i(\tau)B_i\right)y(\tau),
		\qquad y(0)=y_0\in Z,
	\end{align*}
	satisfies $y(\tau)\in Z$ for every $\tau\in[0,T]$. 
 In addition, if
	\begin{align}\label{eq:rank-certificate}
		\mathscr K(M;B_1,\ldots,B_s)f=V
		\qquad\text{for every }f\in Z\setminus\{0\},
	\end{align}
	then $e^{tM}f\in\intt_V Z$ for every $t>0$
	and every $f\in Z\setminus\{0\}$.
\end{proposition}

\begin{proof}
	Fix $f\in Z\setminus\{0\}$ and $t>0$, and set
	$W_i(r)=e^{-rM}B_i e^{rM}f.$
	Suppose that a linear functional $\ell:V\to\R$ satisfies
	$\ell(W_i(r))=0$ for every $i$ and $0<r<t$.
	By differentiating this identity and taking $r\to0^+$,
	we obtain
	\begin{align*}
		0=\ell\bigl(W_i^{(k)}(0)\bigr)
		=(-1)^k\ell\bigl(\operatorname{ad}_M^k(B_i)f\bigr),
		\qquad k\ge0.
	\end{align*}
	By \eqref{eq:rank-certificate}, we have $\ell=0$.
	Hence the vectors $W_i(r)$ span $V$.
	Since $e^{tM}$ is invertible, the vectors
	$
		H_i(r):=e^{tM}W_i(r)=e^{(t-r)M}B_i e^{rM}f
	$
	also span $V$.
	
	Choose $i_1,\ldots,i_m$ and $r_1,\ldots,r_m\in(0,t)$
	such that $H_{i_1}(r_1),\ldots,H_{i_m}(r_m)$ form a basis.
	For each $j$, choose a continuous nonnegative function $\psi_j$
	with integral one and support in a small interval about $r_j$.
	By continuity of $H_{i_j}$, the vectors
	$\int_0^t H_{i_j}(r)\psi_j(r)\,dr$
	approach the chosen basis vectors as these intervals shrink.
	Choose the intervals small enough that the integrated vectors
	remain linearly independent.
	
	For $q=(q_1,\ldots,q_m)\in\R^m$, define
	\begin{align*}
		u_i^q(r)=\sum_{j:i_j=i}q_j\psi_j(r).
	\end{align*}
	Let $y_q$ be the solution with these controls and initial value $f$,
	and set $F(q)=y_q(t)$.
	Since
	\begin{align*}
		\max_i\|u_i^q\|_\infty
		\le\sum_{j=1}^m|q_j|\|\psi_j\|_\infty
		\longrightarrow0
		\qquad(q\to0),
	\end{align*}
	the preservation assumption implies $F(q)\in Z$
	for all sufficiently small $q$.
	
	The map $F$ is continuously differentiable, and
	$F(0)=e^{tM}f$.
	Differentiating the equation with respect to $q_j$ at $q=0$,
	we find that $v_j(r):=\left.\partial_{q_j}y_q(r)\right|_{q=0}$
	satisfies
	\begin{align*}
		\dot v_j(r)=Mv_j(r)+\psi_j(r)B_{i_j}e^{rM}f,
		\qquad v_j(0)=0.
	\end{align*}
	Multiplying by $e^{-rM}$ and integrating, we obtain
	\begin{align}\label{eq:abstract-endpoint-differential}
		\frac{\partial F}{\partial q_j}(0)
		=\int_0^t e^{(t-r)M}B_{i_j}e^{rM}f\,\psi_j(r)\,dr.
	\end{align}
	These vectors are linearly independent by our choice of $\psi_j$,
	so $DF(0)$ is invertible.
	By the inverse function theorem, the endpoints $F(q)$
	for sufficiently small $q$ contain an open neighborhood
	of $F(0)=e^{tM}f$.
	Since these endpoints belong to $Z$, we conclude that
	$e^{tM}f\in\intt_V Z$.
\end{proof}

\begin{proposition}\label{prop:spectral-rank}
Suppose $M$ is symmetric and has an orthonormal eigenbasis
$$v_0=p,v_1,\ldots,v_{m-1},$$ with simple eigenvalues
$\lambda_0>\lambda_1>\cdots>\lambda_{m-1}$.
Suppose all nonzero ordered differences $\lambda_a-\lambda_b$ are
distinct, $\sum_i B_i=c\Id$ with $c\ne0$, and for every $k\ge1$
there is $i$ with $\ip{v_k}{B_ip}\ne0$.
Then
\[
 \mathscr K(M;B_1,\ldots,B_s)f=V
                     \quad\text{whenever }\ip pf\ne0.
\]
\end{proposition}
\begin{proof}
For vectors $u,v$, the notation $uv^*$ means the rank-one map
$z\mapsto\ip vz u$.
The operators $v_av_b^*$ form a basis of $\operatorname{End}(V)$ and
are eigenvectors of $\operatorname{ad}_M$ with eigenvalue
$\lambda_a-\lambda_b$.
Fix $k\ge1$.  On the finite set consisting of zero and the nonzero
ordered differences, polynomial interpolation gives a polynomial
$q_k$ equal to one at $\lambda_k-\lambda_0$ and zero at all other
values.  It can be chosen of degree at most $m(m-1)\le m^2-1$.
Distinctness of the nonzero differences implies
\[
 q_k(\operatorname{ad}_M)B_i=\ip{v_k}{B_ip}\,v_kp^*.
\]
For an index $i$ with nonzero coefficient this proves
$v_kp^*\in\mathscr K(M;B_1,\ldots,B_s)$.
Its value at $f$ is the nonzero multiple $\ip pf v_k$.
Also $\Id\in\mathscr K$ because $\sum_iB_i=c\Id$.
Thus $\mathscr Kf$ contains $v_1,\ldots,v_{m-1}$ and $f$;
the nonzero $p$ coordinate of $f$ makes these vectors span $V$.
\end{proof}

\begin{proposition}\label{prop:ma-interior}
	Let $\mathsf M$ be defined by \eqref{eq:exchange-matrix},
	with parameters $h_i$ and $w_{ij}$ chosen so that all the
	conditions of Lemma~\ref{lem:matrix} hold.
	Then, for every $f\in\mathscr C_{n,d}\setminus\{0\}$
	and every $t>0$, we have
	\begin{align*}
		e^{t\mathsf M}f
		\in\intt_{V_{n,d}^{\ma}}\mathscr C_{n,d}.
	\end{align*}
\end{proposition}

\begin{proof}
The identities $\sum_iN_i=d\Id$ and $d\ge1$, together with
Lemma~\ref{lem:matrix}, give all the hypotheses of
Proposition~\ref{prop:spectral-rank} with $B_i=N_i$ and $M=\mathsf M$.
For nonzero $f\in\mathscr C_{n,d}$ the coefficient vector is
nonnegative and the Perron vector is strictly positive, so
$\ip pf>0$.
The spectral certificate therefore gives
$\mathscr K(\mathsf M;N_1,\ldots,N_n)f=V_{n,d}^{\ma}$.
Lemma~\ref{lem:controlled-preservers} verifies the preservation
assumption of Proposition~\ref{prop:rank-regularization}, with
controls of both signs.  That proposition proves the required
interiority for every positive time.
\end{proof}

\subsection{Normalization and the matrix criterion}
Fix the matrix and the positive unit vector $p$ in
Lemma~\ref{lem:matrix}. Let $\lambda_0$ be its largest eigenvalue, and set
\begin{equation}\label{eq:ma-perron-data}
A=\mathsf M-\lambda_0\Id,\qquad
\eta(f)=\langle p,f\rangle,\qquad
\mathcal A_p=\{f\in V_{n,d}^{\ma}:\eta(f)=1\},\qquad
C_p=\mathscr C_{n,d}\cap\mathcal A_p.
\end{equation}
The functional $\eta$ is positive on every nonzero element of
$\mathscr C_{n,d}$, because every coordinate of $p$ is positive.

\begin{lemma}\label{lem:ma-normalized-flow}
	The set $C_p$ is nonempty, closed, and bounded in $\mathcal A_p$.
	Every $f\in\mathcal A_p$ has a unique expression
$f=p+\sum_{k=1}^{m-1}z_kv_k.$ 
	In these coordinates, $e^{tA}$ acts by $z\mapsto e^{-tD}z$,
	where $D$ is the positive diagonal matrix
	\begin{align}\label{eq:ma-diagonal}
		D=\operatorname{diag}(\lambda_0-\lambda_1,\ldots,
		\lambda_0-\lambda_{m-1}).
	\end{align}
	Moreover, for every $t>0$, we have
$e^{tA}C_p\subseteq\operatorname{int}_{\mathcal A_p}C_p.$
\end{lemma}

\begin{proof}
	Write $p=\sum_I p_Ie_I$, where every $p_I>0$.
	Since $e_I\in\mathscr C_{n,d}$ and $\eta(e_I)=p_I$,
	we have $e_I/p_I\in C_p$. Thus $C_p$ is nonempty.
	Since $\mathscr C_{n,d}$ is closed,
	$C_p=\mathscr C_{n,d}\cap\mathcal A_p$ is closed in $\mathcal A_p$.
	For $f=\sum_I a_Ie_I\in C_p$, we have
	\begin{align*}
		a_I\ge0,\qquad \sum_I p_Ia_I=1.
	\end{align*}
	Hence $0\le a_I\le1/p_I$ for every $I$, so $C_p$ is bounded.
	
	Since $p,v_1,\ldots,v_{m-1}$ form an orthonormal basis
	and $\langle p,f\rangle=1$ for $f\in\mathcal A_p$,
	every such $f$ has the stated unique expression.
	By $A=\mathsf M-\lambda_0\Id$, we have
	$Ap=0$ and $Av_k=(\lambda_k-\lambda_0)v_k$. Therefore,
	\begin{align*}
		e^{tA}f
		=p+\sum_{k=1}^{m-1}e^{-t(\lambda_0-\lambda_k)}z_kv_k.
	\end{align*}
	Thus $e^{tA}$ preserves $\mathcal A_p$ and acts in these
	coordinates by $z\mapsto e^{-tD}z$.
	Since $\lambda_0>\lambda_k$ for $1\le k\le m-1$,
	all diagonal entries of $D$ are positive.
	
	Finally, let $f\in C_p$ and $t>0$.
	Since $\eta(f)=1$, we have $f\ne0$.
	By Proposition~\ref{prop:ma-interior},
	$e^{t\mathsf M}f$ lies in the interior of $\mathscr C_{n,d}$
	in $V_{n,d}^{\ma}$.
	Multiplication by a positive scalar is a homeomorphism
	mapping this cone onto itself. Hence
	\begin{align*}
		e^{tA}f=e^{-t\lambda_0}e^{t\mathsf M}f
		\in\operatorname{int}_{V_{n,d}^{\ma}}\mathscr C_{n,d}.
	\end{align*}
	We also have $e^{tA}f\in\mathcal A_p$.
	Intersecting an open neighborhood of this point contained
	in $\mathscr C_{n,d}$ with $\mathcal A_p$, we obtain a
	relative open neighborhood contained in $C_p$.
	Therefore $e^{tA}f\in\operatorname{int}_{\mathcal A_p}C_p$.
\end{proof}

\begin{proof}[Proof of Theorem~\ref{thm:ma-main}]
	The cases $d=0$ and $d=n$ were proved at the beginning of
	this section. Assume $1\le d\le n-1$ and set $m=\binom nd\ge2$.
	
	By Lemma~\ref{lem:ma-normalized-flow}, the coordinate map
	$\Psi_p:\mathcal A_p\to\mathbb R^{m-1}$ defined by
	\begin{align*}
		\Psi_p\left(p+\sum_{k=1}^{m-1}z_kv_k\right)
		=(z_1,\ldots,z_{m-1})
	\end{align*}
	is an affine isomorphism.
	Set $C=\Psi_p(C_p)$ and $Y=\mathbb R^{m-1}$.
	The same lemma shows that $C$ is nonempty, closed, and bounded,
	that $D$ in \eqref{eq:ma-diagonal} has positive diagonal entries,
	and that
	\begin{align*}
		e^{-tD}C
		=\Psi_p(e^{tA}C_p)
		\subseteq\operatorname{int}C
		\qquad(t>0).
	\end{align*}
	
	We apply Theorem~\ref{thm:main-criterion} with these choices
	and with $B=0$ and $\mathbf c=0$.
	Its differential equation reduces to $\dot z=-Dz$,
	whose solution is $z(t)=e^{-tD}z(0)$.
	Thus the displayed inclusion verifies condition~(i).
	Condition~(ii) holds because
	$\operatorname{int}Y=Y=\mathbb R^{m-1}$.
	The theorem therefore implies that $C$ is homeomorphic
	to $\overline B^{m-1}$.
	Since $\Psi_p$ is an affine isomorphism, the same holds for $C_p$.
	
	It remains to relate $C_p$ to $\mathcal M_{n,d}$.
	For every nonzero $h=\sum_I a_Ie_I\in\mathscr C_{n,d}$,
	the coefficients $a_I$ are nonnegative and not all zero.
	Since every $p_I>0$, we have
	\begin{align*}
		h(\mathbf1)=\sum_I a_I>0,
		\qquad
		\eta(h)=\sum_I p_Ia_I>0.
	\end{align*}
	Define $F:\mathcal M_{n,d}\to C_p$ and
	$G:C_p\to\mathcal M_{n,d}$ by
	\begin{align*}
		F(f)=\frac{f}{\eta(f)},
		\qquad
		G(g)=\frac{g}{g(\mathbf1)}.
	\end{align*}
	Positive rescaling preserves $\mathscr C_{n,d}$, and
	$\eta(F(f))=1$ and $G(g)(\mathbf1)=1$.
	Thus both maps are well defined.
	They are continuous because their denominators are positive
	continuous functions.
	Using $f(\mathbf1)=1$ and $\eta(g)=1$, we obtain
	\begin{align*}
		G(F(f))=\frac{f}{f(\mathbf1)}=f,
		\qquad
		F(G(g))=\frac{g}{\eta(g)}=g.
	\end{align*}
	Hence $F$ and $G$ are mutually inverse homeomorphisms.
	Therefore
	\begin{align*}
		\mathcal M_{n,d}\cong C_p\cong\overline B^{m-1},
		\qquad m-1=\binom nd-1,
	\end{align*}
	as required.
\end{proof}

\appendix
\section{The spectral parameter construction}\label{app:matrix}

We prove Lemma~\ref{lem:matrix}.  All eigenvectors below are vectors of
coefficients in the monomial basis $e_I=x^I$, $I\in\binom{[n]}d$.
The proof separates three tasks: obtaining simple eigenvalues, excluding
identities between eigenvalue differences, and detecting each non-Perron
eigendirection by an occupation operator.

\subsection{Analytic branches and their Taylor coefficients}
At $w_{ij}=0$, the matrix is diagonal, with entries
\begin{equation}\label{eq:diagonal-energies}
 a_I=\sum_{i\in I}h_i.
\end{equation}
For distinct $I,J$, $a_I-a_J$ is a nonzero linear polynomial in $h$.
Avoiding the finitely many hyperplanes where one of these differences
vanishes makes every eigenvalue simple.
Fix a $d$-element set $P$ and restrict the fields to the nonempty open set
\begin{equation}\label{eq:field-chamber}
 \mathcal U_P=\{h:h_i>h_j\text{ whenever }i\in P,\ j\notin P\}.
\end{equation}
Then $a_P>a_I$ for $I\ne P$: the equally many elements lost from $P$ each
have field larger than any element gained outside it.  Removing the
finitely many hyperplanes leaves a nonempty open subset of this chamber.

We use the standard analytic perturbation fact that an isolated simple
eigenvalue of a real-analytic symmetric matrix family has a local
real-analytic eigenvalue branch and a normalized real-analytic eigenvector
branch; see \cite[Chapter II]{Kato}.  A function is real-analytic if it is
locally equal to a convergent power series in its parameters.  For the
simple-eigenvalue case this fact also follows immediately from the analytic
implicit function theorem applied to the eigenvalue equation and a
normalization.  Thus, near $w=0$, write
\begin{equation}\label{eq:eigen-branches}
 \lambda_I(h,w),\qquad v_I(h,w),\qquad
 \lambda_I(h,0)=a_I,\quad v_I(h,0)=e_I,
\end{equation}
with $\norm{v_I}=1$ and its $I$ coordinate positive.  For sufficiently
small real $w$, all branches stay simple and the $P$ branch stays largest.

We need the following more explicit consequence, including at special
fields where eigenvalues other than the two being followed may coincide.
\begin{lemma}\label{lem:rational-coefficients}
The coefficients of the Taylor expansions in $w$ of $\lambda_I$ and $v_I$
are rational functions of the fields $h_i$.  Their denominators involve
only products of the differences $a_I-a_J$, $J\ne I$.
For the scalar $\ip{v_P}{A v_I}$, where $A$ is a fixed diagonal matrix,
the denominators involve only such differences for the branches $P$ and
$I$.  Consequently these coefficients can be specialized whenever those
two diagonal eigenvalues are individually simple, regardless of collisions
among the other diagonal eigenvalues.
\end{lemma}
\begin{proof}
Write the matrix as $D(h)+W(w)$, where $D(h)$ is diagonal and $W$ is linear
in $w$ with zero diagonal.  First use the eigenvector normalization
$u_I(I)=1$ rather than unit length, and write $\lambda_I=a_I+\delta_I$.
The eigenvalue equation reads
\[
 \delta_I=\sum_{J\ne I}W_{IJ}u_I(J),
 \qquad
 (a_I-a_J)u_I(J)=(Wu_I)(J)-\delta_I u_I(J)\quad(J\ne I).
\]
Comparing homogeneous terms in the edge parameters recursively determines
every coefficient, beginning with $u_I=e_I$ and $\delta_I=0$ at $w=0$.
At each recursion step the only divisions are by $a_I-a_J$.
Passing to unit normalization divides by
$\sqrt{1+\sum_{J\ne I}u_I(J)^2}$, whose power series at $w=0$ has rational
numerical coefficients and constant term $1$.  This introduces no further
field denominators.  Multiplication of the two eigenvector series proves
the last assertion.  The analytic implicit function theorem applies to
each selected simple branch even when the remaining eigenvalues collide,
so evaluation of these coefficients at such fields is legitimate.
\end{proof}

\subsection{No unavoidable coincidences between eigenvalue differences}
For a subset $I$, let $\one_I\in\R^n$ be its incidence vector.
Take distinct ordered pairs $(S,T)\ne(U,V)$, with $S\ne T$ and $U\ne V$,
and consider
\begin{equation}\label{eq:gap-function}
 F_{S,T,U,V}(h,w)=\lambda_S-\lambda_T-\lambda_U+\lambda_V.
\end{equation}
We show that some Taylor coefficient in $w$ is a nonzero rational function
of $h$.  If $\one_S-\one_T\ne\one_U-\one_V$, its constant coefficient
$a_S-a_T-a_U+a_V$ is already a nonzero linear polynomial.

In the remaining case, the positive and negative coordinates of the common
incidence difference give disjoint nonempty sets $A,B$, and there are sets
$C,C'$ disjoint from $A\cup B$ such that
\[
 S=A\sqcup C,\quad T=B\sqcup C,\qquad
 U=A\sqcup C',\quad V=B\sqcup C'.
\]
Here $|A|=|B|$, $|C|=|C'|$, and $C\ne C'$ because the ordered pairs are
distinct.  Choose $a\in A$, $b\in B$, and $c\in C\setminus C'$.
We compute the coefficient of $w_{ab}w_{ac}w_{bc}$ by setting all other
edge variables to zero.  For this computation take fields with all the
required diagonal differences nonzero, and set
\[
 F_a=\frac{2}{(h_a-h_b)(h_a-h_c)},\quad
 F_b=\frac{2}{(h_b-h_a)(h_b-h_c)},\quad
 F_c=\frac{2}{(h_c-h_a)(h_c-h_b)}.
\]
If a basis state occupies just $a$ in the triangle $\{a,b,c\}$, the
coefficient in its eigenvalue branch is $F_a$.
To verify this, the recursion in Lemma~\ref{lem:rational-coefficients}
gives the third-degree eigenvalue term
\[
 \sum_{J,L\ne I}
 \frac{W_{IJ}W_{JL}W_{LI}}{(a_I-a_J)(a_I-a_L)}.
\]
The monomial in question occurs on the two oriented walks once around
the triangle.  Their denominator products are both
$(h_a-h_b)(h_a-h_c)$, giving the factor $2$.
If a state occupies two vertices and is missing $a$, the same coefficient
is $F_a$: the two energy differences each change sign, leaving their
product unchanged.  The analogous statements hold for $b$ and $c$.

The occupancies of $S,T,U,V$ in this triangle are, respectively,
$\{a,c\},\{b,c\},\{a\},\{b\}$.  Thus the coefficient in
\eqref{eq:gap-function} is
\[
 F_b-F_a-F_a+F_b=2(F_b-F_a),
 \qquad
 F_b-F_a=
 -\frac{2(h_a+h_b-2h_c)}{(h_a-h_b)(h_a-h_c)(h_b-h_c)}.
\]
This is a nonzero rational function.  Every nontrivial equality between
ordered gaps is therefore excluded by avoiding the zero set of a
nonidentically zero analytic function, once suitable fields are fixed.

\subsection{A nonvanishing coupling coefficient}
Fix $S\ne P$, and put
\[
 A=P\setminus S,\qquad B=S\setminus P,\qquad C=P\cap S,
 \qquad |A|=|B|=k\ge1.
\]
Let $N_A=\sum_{i\in A}N_i$.  We will show that
\begin{equation}\label{eq:coupling-function}
 g_{P,S}(h,w)=\ip{v_P(h,w)}{N_Av_S(h,w)}
\end{equation}
has a Taylor coefficient in the edge variables that is a nonzero rational
function of $h$.

Choose numbers $h>0$ and $L>kh$, and specialize the fields to
\begin{equation}\label{eq:special-fields}
 h_i=\begin{cases}
 h,&i\in A,\\
 0,&i\in B,\\
 L,&i\in C,\\
 -L,&i\notin P\cup S.
 \end{cases}
\end{equation}
At zero exchange, $a_S=|C|L$ and $a_P=|C|L+kh$ are both simple in the
\emph{full} matrix.  To check this, for any $d$-set $I$ write
\[
 a_I-|C|L=-j_0L+j_1h,
 \quad j_0=|C|-|I\cap C|+|I\setminus(P\cup S)|\ge0,
 \quad 0\le j_1=|I\cap A|\le k.
\]
If $j_0\ge1$, this expression is negative.  If $j_0=0$, then $I$ contains
all of $C$ and no point outside $P\cup S$; the expressions $0$ and $kh$
occur only for $I=S$ and $I=P$, respectively.  Hence the two endpoint
eigenvalues are isolated from every other diagonal eigenvalue.  This
argument permits degeneracies among the intermediate eigenvalues.

Keep only the following edge parameters: give each edge inside $A$ and
inside $B$ the weight $q/2$, and each edge between $A$ and $B$ the weight
$b$; set all other weights to zero.  These parameters are used only to
certify a nonzero Taylor coefficient.  Strict positivity of every edge
weight will be arranged later.
Choose $q,b$ in a sufficiently small neighborhood of zero that the two
endpoint eigenvalues remain individually simple in the full matrix.
Such a neighborhood exists by their positive separation just proved and
continuity of the eigenvalues.  We may also assume $q>0$ and $h>2kq$
when computing first as a power series in $b$.

The subspace with all of $C$ occupied, no outside point occupied, and exactly
$k$ occupied points in $A\cup B$ is invariant.  Within it, the subspace
invariant under all permutations of $A$ and all permutations of $B$ has
the following orthonormal basis:
\begin{equation}\label{eq:uniform-basis}
 u_j=\binom kj^{-1}
 \sum_{\substack{A'\subseteq A,\ |A'|=j\\B'\subseteq B,\ |B'|=k-j}}
       e_{C\cup A'\cup B'},\qquad 0\le j\le k.
\end{equation}
There are $\binom kj^2$ terms in this sum, explaining the normalization.
The branches from $e_S=u_0$ and $e_P=u_k$ stay in this invariant subspace:
the matrix commutes with those permutations and a simple normalized
eigenvector continuing either endpoint is fixed by them.  Equivalently,
the uniqueness of the eigenvector branch and the invariant block containing
the initial vector give the same conclusion.

After subtracting the constant $|C|L$ from the diagonal, the restriction
is the symmetric tridiagonal matrix $J$ with
\begin{equation}\label{eq:Jacobi}
 J_{jj}=E_j=hj+qj(k-j),\qquad
 J_{j,j+1}=J_{j+1,j}=b(j+1)(k-j).
\end{equation}
For the diagonal contribution of the exchanges, a state in layer $j$ has
$j(k-j)$ possible exchanges in $A$ and the same number in $B$, each of
weight $q/2$.  Thus the uniform vector has eigenvalue $qj(k-j)$ under
these within-layer exchanges.
For the off-diagonal entry, each state in layer $j$ has $(k-j)^2$ exchanges
to layer $j+1$, and each state in layer $j+1$ has $(j+1)^2$ exchanges
back.  The normalized entry is the square root of their product, times
$b$, giving \eqref{eq:Jacobi}.
On this subspace $N_A=\diag(0,1,\ldots,k)$.

Let $v^{(0)}(b)$ and $v^{(k)}(b)$ be the unit eigenvector branches of $J$
from $u_0$ and $u_k$, with the endpoint coordinate positive.
For fixed small $q$, the leading term of the $j$th coordinate of
$v^{(0)}$ is
\[
 b^j\frac{\prod_{r=0}^{j-1}(r+1)(k-r)}
             {\prod_{\ell=1}^{j}(E_0-E_\ell)},
\]
and that of $v^{(k)}$ is
\[
 b^{k-j}\frac{\prod_{r=j}^{k-1}(r+1)(k-r)}
                   {\prod_{\ell=j}^{k-1}(E_k-E_\ell)}.
\]
Empty products are $1$.  These formulas follow by recursively comparing
powers of $b$ in
$(\lambda-E_j)v_j=bj(k-j+1)v_{j-1}+b(j+1)(k-j)v_{j+1}$.
Reaching layer $j$ from an endpoint first requires exactly the distance
from that endpoint in off-diagonal steps; at that first possible order
there is only the monotone path.  The eigenvalue shift is $O(b^2)$ and
normalizing to unit length changes an endpoint-normalized vector by
$1+O(b^2)$, so neither affects these leading terms.
Here $O(b^a)$ means a term bounded in absolute value by a constant times
$|b|^a$ as $b\to0$, with the other parameters fixed in the stated
neighborhood.

Multiplying the two coordinate expansions and using the diagonal entries
of $N_A$ gives
\begin{equation}\label{eq:coupling-leading}
 \ip{v^{(k)}}{N_Av^{(0)}}
 =(k!)^2b^k\sum_{j=0}^k
 \frac{j}{\displaystyle
       \prod_{\ell=j}^{k-1}(E_k-E_\ell)
       \prod_{\ell=1}^{j}(E_0-E_\ell)}
 +O(b^{k+1}).
\end{equation}
The common numerator $(k!)^2$ is
$\prod_{r=0}^{k-1}(r+1)(k-r)$.

\begin{lemma}\label{lem:finite-difference}
For $E_j=hj+qj(k-j)$ and $h>2kq>0$,
\begin{equation}\label{eq:sum-identity}
 \sum_{j=0}^k
 \frac{j}{\displaystyle
       \prod_{\ell=j}^{k-1}(E_k-E_\ell)
       \prod_{\ell=1}^{j}(E_0-E_\ell)}
 =-\frac{q^{k-1}}{h\prod_{r=1}^{k-1}(h^2-r^2q^2)}.
\end{equation}
\end{lemma}
\begin{proof}
The differences are
$E_k-E_\ell=(k-\ell)(h-q\ell)$ and
$E_0-E_\ell=-\ell(h+q(k-\ell))$.
Put $a=h/q$.  On bringing the sum to the common denominator
$q^k a\prod_{r=1}^{k-1}(a^2-r^2)$, its numerator becomes
\[
 \sum_{j=0}^k\frac{(-1)^j}{j!(k-j)!}
       j\prod_{\ell=1}^{k-1}(a-j+\ell).
\]
Indeed, the factors left after cancellation from the common denominator
are exactly $\prod_{\ell=1}^{k-1}(a-j+\ell)$.
The polynomial $F(j)=j\prod_{\ell=1}^{k-1}(a-j+\ell)$ has degree $k$
and leading coefficient $(-1)^{k-1}$.
For the forward difference $\Delta F(j)=F(j+1)-F(j)$, induction gives
$\Delta^kF(0)=\sum_{j=0}^k(-1)^{k-j}\binom kj F(j)$.
Successive differences lower the degree by one and multiply the leading
coefficient by the current degree, so
$\Delta^kF(0)=k!(-1)^{k-1}$.
The displayed numerator is consequently $-1$.
Converting the denominator back from $a$ to $h,q$ proves
\eqref{eq:sum-identity}, also for $k=1$ with the indicated empty-product
convention.
\end{proof}

Equations~\eqref{eq:coupling-leading}--\eqref{eq:sum-identity} show that the
specialized coupling is
\begin{equation}\label{eq:coupling-nonzero}
 -\frac{(k!)^2q^{k-1}b^k}{h\prod_{r=1}^{k-1}(h^2-r^2q^2)}+O(b^{k+1}).
\end{equation}
In particular, writing $[b^uq^v]$ for a coefficient of the joint Taylor
series at $(b,q)=(0,0)$, we obtain
\begin{equation}\label{eq:joint-coefficient}
 [b^kq^{k-1}]\ip{v^{(k)}}{N_Av^{(0)}}
 =-\frac{(k!)^2}{h^{2k-1}}\ne0.
\end{equation}
The equality follows by expanding the rational coefficient of $b^k$ in
\eqref{eq:coupling-nonzero} at $q=0$.  That coefficient is analytic there
because $h>0$; terms of order higher than $b^k$ cannot contribute.

The specialized edge variables are linear functions of $b,q$.
Hence \eqref{eq:joint-coefficient} is a finite linear combination of
Taylor coefficients in the original edge variables of
\eqref{eq:coupling-function}, evaluated at \eqref{eq:special-fields}.
By Lemma~\ref{lem:rational-coefficients} these are rational functions of
the fields, and the specialization is valid because $a_P$ and $a_S$ are
simple in the full diagonal matrix.  The nonzero value proves that at
least one of those rational functions is not identically zero.
This establishes the required nonvanishing certificate for every $S\ne P$.

\subsection{Simultaneous parameter selection and completion}
There are only finitely many ordered-pair gap conditions
\eqref{eq:gap-function}, and finitely many couplings
\eqref{eq:coupling-function}.  For each gap choose one of the nonzero
rational Taylor coefficients established above, including the constant
coefficient when appropriate.  For each coupling choose a nonzero
rational Taylor coefficient supplied by \eqref{eq:joint-coefficient}.
Select $h\in\mathcal U_P$ for which all subset sums are distinct and all
these chosen coefficients are defined and nonzero.
Such a choice exists: there are finitely many nonzero polynomial numerators
and denominators to avoid; the zero set of a nonzero real polynomial has
empty interior, so their finite union cannot exhaust the nonempty open
chamber.  This also explains why the special fields used for one coupling
need not simultaneously serve all other couplings.

Fix these fields.  On a sufficiently small connected ball of edge parameters
around $w=0$, all eigenvalues are simple, their analytic branches are
defined, and the $P$ branch is largest.
Every selected gap function and coupling function is nonidentically zero
on this ball, because its chosen Taylor coefficient is nonzero.
A nonidentically zero real-analytic function on a connected open set
cannot vanish on a nonempty open subset; this is the real-analytic identity
principle, obtained by continuation of its convergent local power series.
Its zero set is closed with empty interior.  A finite union of these zero
sets cannot cover the nonempty open part of the ball where every $w_{ij}>0$.
Choose $w$ in that positive part avoiding them all.

The resulting symmetric matrix has simple eigenvalues, and all
its nonzero ordered eigenvalue differences are distinct.
Its off-diagonal graph is the connected exchange graph.
By the argument following Theorem~\ref{thm:perron},
a unit eigenvector for its largest eigenvalue can be chosen
strictly positive.
The $P$ branch corresponds to this eigenvalue.
Since its chosen unit eigenvector has a positive $P$ coordinate,
its sign already agrees with that of the positive Perron vector.
For every other branch $S$,
\[
 0\ne\ip{v_P}{N_{P\setminus S}v_S}
   =\sum_{i\in P\setminus S}\ip{v_S}{N_iv_P},
\]
where we used the symmetry of $N_i$.  At least one summand is nonzero.
Relabel the branches in decreasing eigenvalue order, with $v_0=p=v_P$.
All three requirements of Lemma~\ref{lem:matrix} now hold, completing
its proof.

\section*{Declaration of AI usage}
During the preparation of this manuscript, the authors used generative AI tools to assist in exploring possible approaches, checking technical details, and improving the exposition. The authors verified the mathematical arguments, computations, proofs, and results and take full responsibility for the content of this manuscript.

\end{document}